\documentclass[12pt]{amsart}
\usepackage{amsfonts}

\usepackage[all,cmtip]{xy}
\usepackage{amsmath}
\usepackage{amsthm}
\usepackage{amssymb}
\usepackage{enumitem}

\newtheorem*{definition*}         {Definition}

\theoremstyle{remark}

\newcommand*{\Q}{\mathbb{Q}}

\newcommand*{\Z}{\mathbb{Z}}
\newcommand*{\N}{\mathbb{N}}

\newcommand*{\GG}{\mathbb{G}}
\newcommand*{\R}{\mathbb{R}}
\newcommand*{\C}{\mathbb{C}}

\newcommand*{\PP}{\mathbb{P}}

\usepackage{amscd,amssymb,amsmath}
\usepackage{color}
\definecolor{purple}{rgb}{0.59, 0.44, 0.84}

\def\hash{\#}

\begin{document}


%
%
\newcount\refCount
\def\newref#1 {\advance\refCount by 1
\expandafter\edef\csname#1\endcsname{\the\refCount}}
\newref ANDRE
\newref AWH
\newref AX
\newref BALSER
\newref BATEMAN
\newref BEURLING
\newref BERGER
\newref BKPLMM
\newref BNZ
\newref BCFNDAS
\newref DLMF
\newref DRIESSPEISSEGGER
\newref DYNKIN
\newref FUGLEDE
\newref FUGLEDEBALL
\newref GELFOND
\newref IK
\newref JK
\newref KIRBY
\newref KOLCHIN
\newref LANG
\newref MANGOUBIROSENBLATT
\newref MANTOVA
\newref MARKER
\newref OLVER
\newref PW
\newref PZ
\newref PUTSINGER
\newref RAYNAUD
\newref SCHNEIDER
\newref SHIDLOVSKII
\newref SIEGEL
\newref SIEGELBOOK
\newref WALDSCHMIDT
\newref WASOW
\newref ZAKHAROV

\newcount\derandomCount
\def\derandom#1 {\advance\derandomCount by 1
\expandafter\edef\csname#1\endcsname{\the\derandomCount}}

\title{Counting Bessel zeros and applications}
\author{Gal Binyamini, Avner Kiro, 
Gady Kozma, and Jonathan Pila}

\maketitle

{\bf Abstract.\/} We establish o-minimal definability 
results and develop point-counting machinery applicable
to the zeros of a Bessel function of fixed non-negative 
real parameter.
We also give an Ax-Schanuel theorem for Bessel functions 
of fixed (complex) parameters and exponential functions. 
We combine these to prove a sparsity result 
(weak Manin-Mumford) for points on an algebraic 
variety whose coordinates are zeros of a Bessel 
function of fixed (non-negative real) parameter. We apply 
this result to improve density results on spectral sets 
for the unit ball motivated by Fuglede's conjecture.

\bigskip

\section{Introduction}

\bigskip

Let $J_\kappa(x)$ denote the Bessel 
function of the first kind of order 
$\kappa$ (see \S2). In general 
$\kappa\in \C$ but we will mostly be
concerned with $\kappa\in \R_{\ge 0}$ 
and even with the orders
$\kappa=m/2\in \N \cup (\N+1/2)$ 
($\N=\{0,1,2,\ldots\}$) connected with the
spectral theory of the unit ball in $\R^m$.

For $\kappa\in\R$, $J_\kappa(x)$ has 
infinitely many real zeros, symmetrically 
located about $x=0$, and only finitely 
many non-real zeros (see [\DLMF, 10.21(i)]).
We show that, for non-negative
real $\kappa\ne 1/2$, 
$n$-tuples of zeros of $J_\kappa(x)$ 
on an algebraic variety 
$V\subset (\C^\times)^n$
must be quite sparse, 
subject only to the obvious
necessary obstruction: the presence 
of positive-dimensional 
``Bessel special subvarieties'' 
(see Definition 1.1) contained in $V$.
For $\kappa= 1/2$, the zeros are 
precisely $\pi\Z$ (for $\kappa=-1/2$
they are $\pi(\Z+1/2)$)
and our approach breaks down
(the failure can be seen in Lemma 5.2).

Though $J_\kappa$ is sometimes 
well-defined at $x=0$, the variety 
$\C^\times$ comprising the regular 
points of the differential 
equation is the natural ambient.

\medbreak
\noindent
{\bf 1.1. Definition.\/} 
Fix $\kappa\in\C$.
{\it A (Bessel) special point\/} 
is a zero $\beta\in\C^\times$ of 
$J_\kappa$. 
A {\it (Bessel) special point\/} in 
$(\C^\times)^n$ is an $n$-tuple 
$(\beta_1,\ldots, \beta_n)$
of Bessel special points. 
A {\it (Bessel) special subvariety\/} 
in $(\C^\times)^n$ is a non-empty 
subvariety defined by any number of 
equations of the forms: $x_i=\beta_i$, 
where $\beta_i$ is a Bessel special 
point, or $x_j=\pm x_k$ (where $j\ne k$).

\medbreak
Fixing $\kappa\in [0,\infty)$, the 
(positive) special points have asymptotic 
expansions due to McMahon (see e.g. 
[\DLMF, \S10.21(vi)]) which give that, 
as $m\rightarrow+\infty$,
$$
\beta_m=\beta_{\kappa, m}=
\pi\big(m+\kappa/2-1/4\big)+O_\kappa\big(m^{-1}\big),
$$
where $\beta_{\kappa, m}$ is the $m$th positive zero 
(which are all simple) of $J_\kappa$, 
and $m=1, 2,\ldots$. 
We set the {\it height\/}
of $\beta_m$ to be 
$H(\beta_m)=|\beta_m|$,
which thus is roughly linear in $m$ 
(for fixed $\kappa$).
If $T$ is a special subvariety we define
its {\it height\/}, denoted $H(T)$, 
to be the maximum of 1 and the height
of its constant (special point)
coordinates.
A special point is then the same as a 
special subvariety of dimension zero, 
while a special subvariety of
positive dimension $k$ contains ``many''
special points in a height density sense
(asymptotically $2^kH^k/\pi^k$ special 
points up to height $H$). 

By a {\it 
maximal\/} special subvariety $T\subset V$,
for a variety $V$, we mean one maximal for 
inclusion among special subvarieties 
contained in $V$.

\medbreak
\noindent
{\bf 1.2. Theorem.\/} 
{\it Let 
$\kappa\in[0, \infty)]\backslash\{1/2\}$, 
$\epsilon>0$, and let
$V\subset (\C^\times)^n$ be an 
algebraic subvariety. 
Then, for $H\ge 1$, $V$ contains at most
$$
C(V, \kappa, \epsilon)H^\epsilon
$$ 
maximal Bessel special subvarieties $T$ 
with height $H(T)\le H$.

For $\kappa\in \N+1/2$ (except 
$\kappa=\pm 1/2$) the estimate may be 
improved  to $$C(V, \kappa)
(\log (1+H))^{c(V, \kappa)}.$$
\/}


In particular, if $V$ contains no 
positive-dimensional special subvarieties 
then these estimates hold for the number 
of special points on $V$ 
up to height $H$. Our argument does 
not render effective constants
though they are uniform
for all $V\subset (\C^\times)^n$ 
of given dimension and degree.
The poly-logarithmic bound for 
$\kappa\in \N+1/2, \kappa\ne 1/2$,
follows from the proof [\BNZ] of 
Wilkie's conjecture on point-counting 
bounds for sets definable using 
restricted Pfaffian functions and 
(unrestricted) exponential. 

Theorem 1.2 is a weak analogue of the
classical Manin-Mumford conjecture 
(theorem of Raynaud [\RAYNAUD]). 
It is weak in asserting only sparsity of the
(maximal) special subvarieties in $V$ rather
than finiteness. This is due to the fact that
the special points are in general not algebraic numbers 
(by Siegel [\SIEGEL], the non-zero
roots are transcendental when $\kappa\in\Q$),
and so the algebraic structure 
that enables deducing finiteness from sparseness 
in the classical context is absent.

The proof of Theorem 1.2 uses point-counting and functional
transcendence. It follows the basic strategy
introduced in [\PZ], but with some key new features,
including the manner in which special points are counted 
by integer points in a definable set.
Definability of the counting set in an o-minimal structure
entails that suitable related functions (essentially the Hankel
functions with their oscillation ``removed'') are so definable.
Showing this uses summability theory of linear
differential equations and quasi-analyticity;
the relevant o-minimal structure is the
structure $\R_{\mathcal{G},\exp}$ of
{\it multi-summable Gevrey series\/} 
(with exponentiation) introduced 
(and proved o-minimal) by van den Dries-Speissegger 
[\DRIESSPEISSEGGER].

Our functional transcendence result is 
an ``Ax-Schanuel" statement for
Bessel functions and exponentials.
It is a generalization 
of the classical Ax-Schanuel theorem
of Ax [\AX, Theorem 3].
A very special case is a result of Siegel
[\SIEGEL, Satz 9], proved in the context of the
above-mentioned transcendence results.
This result was 
reproved by Kolchin [\KOLCHIN], who introduced 
some group-theoretic ideas that are central to our treatment
(see [\JK] for a generalization).
We combine  Kolchin's results and a study of the
formal solutions and their Stokes phenomenon
with the recent Ax-Schanuel technology [\BCFNDAS] of 
Bl\'azquez-Sanz--Casale--Freitag--Nagloo.

\medbreak

\noindent{\bf 1.3. Definition.\/} 
A {\it (Bessel) weakly special 
subvariety\/} of $(\C^\times)^n$
is a subvariety defined by any number of
equations of either of the two forms: 
$x_i=c_i, c_i\in\C^\times$, and 
$x_j=\pm x_k$ (where $j\ne k$).

\medbreak

Let $K$ be a differential field with 
derivation  $D$ and field of constants 
$\C$. For a field
$M$, subfield $L$, and elements $x_1,\ldots, x_n$
of $M$, tuples  
$\overline{y}_1,\ldots, \overline{y}_m$ of elements
of $M$, matrices $A_1,\ldots, A_{\ell}$
of elements of $M$, tuples of matrices
$\overline{A}_1,\ldots, \overline{A}_k$, etc,
we denote by 
$$
F=L(x_1,\ldots, x_n, \overline{y}_1,\ldots, 
\overline{y}_m, A_1,\ldots, A_{\ell}, 
\overline{A}_{1}, \ldots, \overline{A}_{k},{\rm etc})
$$
the field obtained by adjoining to $L$ all the 
constituent elements and entries of the tuples 
and matrices etc, and by ${\rm td}_LF$ we denote 
the transcendence degree of the extension $F/L$.

{\it Bessel's differential equation of order 
$\kappa$\/} (of which $z=J_\kappa(x)$ is a solution)
is
$$
x^2z''+xz'+(x^2-\kappa^2)z=0,
$$
where $'$ denotes differentiation with respect to $x$. 
We call the equation
{\it generic\/} if $\kappa\not\in \Z+1/2$. In that
case, the differential Galois group is well-known
(see e.g. Kolchin [\KOLCHIN])
to be ${\rm SL}_2(\C)$. Thus, if
$$
W=
\begin{pmatrix}z_1 & z_2 \cr z_1' & z_2'
\end{pmatrix}
$$
is a {\it fundamental matrix\/} formed from two
$\C$-linearly independent solutions $z_1, z_2$
then
$
{\rm td}_{\C(x)}\C(x, W)=3.
$
For non-generic $\kappa$, the
solutions to Bessel's differential equation are
elementary (combinations of exponentials $e^{\pm ix}$
and rational functions of $x$; see \S2), and 
${\rm td}_{\C(x)}\C(x, W)=1$.
Here is our Ax-Schanuel result
for Bessel functions and exponentials.
Below, all the functions are assumed contained in $K$.

\medbreak
\noindent
{\bf 1.4. Theorem.} 
{\it Let $x_1,\ldots, x_n \in K\backslash \C$
and, for each $i$, let $W_i$ be a fundamental
matrix of solutions to Bessel's differential
equation (with respect to $x_i$)
of generic order $\kappa_i\in\C$.
Let further $y_1,\ldots, y_m\in K\backslash \C$
and, for each $j$, let $Y_j$ be a 
fundamental matrix of solutions to Bessel's
differential equation (with respect to $y_j$)
of non-generic order $\nu_j$.
Finally, let $e_1,\ldots, e_{\ell}, t_1,\ldots, t_\ell\in K\backslash \C$
with, for each $k$, $Dt_k=De_k/e_k$.
Let $\overline{x}=(x_1,\ldots, x_n), \overline{W}=
(W_1,\ldots, W_n)$ etc.
Then
$$
{\rm td}_{\C}\C\Big(\overline{t}, \overline{e}, 
\overline{y}, \overline{Y}, \overline{x}, 
\overline{W}\Big) \ge \ell+ m+3n+1
$$
unless either

(A) $t_1,\ldots, t_{\ell}, iy_1,\ldots, iy_m$ are
linearly dependent over $\Q$ modulo $\C$, or

(B) there are $i, j\ne i$ such that $x_i=x_j$
or $x_i=-x_j$
(and $\kappa_i\pm \kappa_j\in\Z$).\/}

\medbreak

We apply the above results to study sets 
$\Lambda\subset \R^m$ of {\it orthogonal
exponentials\/} for the unit ball $\Delta_m\subset\R^m$.
See [\FUGLEDE, \FUGLEDEBALL, \IK, \ZAKHAROV]. 
These are sets such that,
for any $\lambda_1,\lambda_2 \in \Lambda$ with
$\lambda_1 \neq \lambda_2$ we have
$$
\int_{\Delta_m} \exp\big( 2\pi i u \cdot 
(\lambda_1-\lambda_2)\big)du =0.
$$

The orthogonality condition is equivalent to
the distance between any two distinct points of
$\Lambda$ being a zero of $J_{m/2}(2\pi x)$.
Equivalently, $\Lambda^*=2\pi \Lambda$ is a {\it Bessel distance set\/} 
in $\R^m$: all positive distances are zeros of $J_{m/2}$.

It is known that any such set is finite
(Fuglede [\FUGLEDE] states finiteness for $m=2$ and 
in [\FUGLEDEBALL] proves finiteness for all $m\ge 2$). 
For given $m$ it is not known if the size 
of such a set is uniformly bounded.
One can construct an abundance of
such sets of size $m+1$, but Fuglede 
[\FUGLEDE, \FUGLEDEBALL] conjectured that such a 
set of size $m+2$ in $\R^m$ is impossible, observing
that the existence of such a set would entail
some algebraic relation among the positive
zeros of $J_{m/2}$, which he viewed as ``highly unlikely''.

Iosevich and Kolountzakis [\IK] give a
uniform linear bound 
for the size of such $\Lambda$ for $m=2$ in terms of
the infimum of non-zero differences
$$
\inf\{|\lambda_1-\lambda_2|: \lambda_1, \lambda_2 \in \Lambda, 
\lambda_1\neq \lambda_2\}
$$
in $\Lambda$, and show 
further that such sets are sparse in the sense that
$$
\hash \Lambda \cap [-R, R]^2 \le cR^{2/3}.
$$
The exponent $2/3$ is improved to $3/5+\epsilon$ 
by Zakharov [\ZAKHAROV].

If we consider the tuple of
distances associated with $m+2$ element subsets
of $\Lambda$ we generate special points
on the distance variety $V_{m+2, \R^m}\subset 
\R^M, M=(m+2)(m+1)/2$ (whose points are the $M$-tuples of distances between $m+2$ points in $\R^m$).
This set is contained in a strict subvariety of 
$\R^{(m+2)(m+1)/2}$ (see [\BERGER, \S9.7]). 
Using 1.2 we prove

\medbreak
\noindent
{\bf 1.5. Theorem.\/}
{\it Let $m\ge 2$ and $\epsilon >0$.
For a set $\Lambda\subset \R^m$ of orthogonal 
exponentials for the unit ball in $\R^m$ we have
$$
\hash \Lambda \cap [-R, R]^m \le 
c(m, \epsilon)R^{\epsilon}.
$$
For odd $m$ the estimate may be 
improved  to $C(m)(\log R)^{c(m)}.$
}

\medbreak

Note that the  constraints afforded by 
Theorem 1.4 are much less restrictive 
for non-generic orders than for
generic ones. It is remarkable that, 
despite this, the same form holds 
for all non-negative real $\kappa \ne 1/2$
for the weakly special subvarieties in 
Theorem 1.2 (cf Lemma 5.1 versus Lemma 5.2),
giving the same form of estimate in 
Theorem 1.5.

The proof of Theorem 1.5 is in \S7.
Our present argument, relying on 
Theorem 1.2, does not render effective 
constants.

In the final section \S8 of the paper 
we return to the question of possible 
algebraic relations among zeros of the 
Bessel functions relevant to Fuglede's 
conjecture, and more generally those 
of rational order.
For non-generic $J_{m+1/2}$, 
the existence of such relations
is governed by Schanuel's conjecture 
([\LANG, p30], Conjecture 8.1)
on the transcendence properties of the 
exponential function. We show that 
Schanuel's conjecture implies the 
algebraic independence of distinct 
positive zeros and thus precludes the 
existence of a set of orthogonal 
exponentials in $\R^m$ of size $m+2$ 
for odd $m\ge 3$.

For even $m=2k$, the algebraic 
independence of positive zeros of $J_k$ would also rule
out the existence of a set of orthogonal
exponentials of size $m+2$ in $\R^m$. However there seems to be
no ``Bessel-Schanuel conjecture'' formulated in the literature.
Guided by our
Bessel-Ax-Schanuel theorem above and the Siegel-Shidlovskii 
theorems [\SIEGEL, \SHIDLOVSKII], we offer a 
formulation in Conjecture 8.4.

\medbreak
\noindent
{\bf Acknowledgments.\/} 
GB was supported by the Marvin V. and Beverly J. Mielke Endowed Fund and the Infosys Member Fund as a member of the IAS. GB and AK were supported by the European Union (ERC, SharpOS, 101087910) and by Israel Science Foundation grant no. 2067/23.

JP is grateful 
to the Weizmann Institute of Science and 
the Institute for Advanced Study for 
their warm hospitality and financial 
support. JP was supported as a Weston 
Visiting Professor during his visit to 
the WIS, and by the Charles Simonyi 
Endowment as a member of the IAS. 

GK was supported by the Israel Science Foundation and the Binational Science Foundation.

\medskip
\noindent
{\bf AI disclosure.\/} GPT 5.6 sol and 6 Astra were used 
to check various details, including
identifying some references, and to 
review an earlier fully human-written
draft of this paper, leading to various 
suggested corrections that were 
human-implemented. In particular it spotted an
error in an earlier version of the proof of Theorem 3.1 for rational
non-integer parameters, and assisted in fixing this
error. Note that this error did not affect our main 
applications for the Fuglede problem which deal with
integer parameters.

\section{Bessel Basics}

Bessel's differential equation (above)
is of confluent hypergeometric type.
It has an irregular singular point at $x=\infty$
and a branch point at $x=0$.

The  {\it Bessel function of the first kind,\/}
$J_\kappa(x)$, may be defined by
$$
J_\kappa(x)=\Big({x\over 2}\Big)^\kappa 
\sum_{j=0}^\infty{(-1)^j\over 
\Gamma(j+1)\Gamma(j+\kappa+1)} 
\Big({x\over 2}\Big)^{2j}.
$$
(Note that $\Gamma(z)^{-1}$ is an entire function so the sum is always an entire function.)

For integer $\kappa$ the function $J_\kappa(x)$ is 
single valued (and even entire for non-negative integer $\kappa$), 
otherwise it is multivalued, due (only) to the factor $(x/2)^\kappa$. 
In particular, the zeros of $J_{\kappa}$ are well-defined.

For non-integer $\kappa$, a second $\C$-linearly independent 
solution is provided
by $J_{-\kappa}(x)$, and one defines the 
{\it Bessel function of the
second kind,\/} $Y_\kappa$, to be 
$$
Y_\kappa(x)=\frac{J_\kappa(x)\cos(\pi\kappa)-J_{-\kappa}(x)}{\sin(\pi\kappa)}.
$$

For integer $\kappa$ the functions $J_\kappa, J_{-\kappa}$
are dependent ($J_\kappa=(-1)^\kappa J_{-\kappa}$) and one defines the
{\it Bessel function of the second kind\/} $Y_\kappa$
by the limit
$$
Y_\kappa(x)=\lim_{\alpha\to\kappa}Y_\alpha(x),\quad \kappa\in \mathbb{Z}.
$$

In all cases, $J_\kappa, Y_\kappa$ are $\C$-linearly independent solutions to 
Bessel's d.e. of order $\kappa$.
The Wronskian determinant is given by
$$
{\rm Wr}(J_\kappa, Y_\kappa)=\det
\begin{pmatrix}
J_\kappa(x) & Y_\kappa(x)  \cr
 J_\kappa '(x) & Y_\kappa'(x)
\end{pmatrix}={2\over \pi x}.
$$

Finally, for all $\kappa$, {\it Hankel functions\/}
(also known as {\it Bessel functions of the third kind\/})
are further canonical solutions of Bessel's differential
equation defined by
$$
H^{(1)}_\kappa(x)= J_\kappa(x)+iY_\kappa(x),
\quad 
H^{(2)}_\kappa(x)= J_\kappa(x)-iY_\kappa(x).
$$
The principal branches have the asymptotic series 
$$
H^{(1)}_\kappa(x)=
\sqrt{2\over \pi x}
e^{i(x-\kappa\pi/2-\pi/4)}\sum_{m=0}^\infty
{({1\over 2}-\kappa)_m
{({1\over 2}+\kappa)_m}
\over
m!(2ix)^m},
$$
$$
H^{(2)}_\kappa(x)=
\sqrt{2\over \pi x}
e^{-i(x-\kappa\pi/2-\pi/4)}\sum_{m=0}^\infty
{({1\over 2}-\kappa)_m
{({1\over 2}+\kappa)_m}
\over
m!(-2ix)^m},
$$
as $|x|\to\infty$ in the sectors $\text{arg}\,x\in(-\pi,2\pi)$  and $\text{arg}\,x\in(-2\pi,\pi)$, 
where the Pochhammer symbol $(x)_m$ denotes the
``raising factorial''
$$
(x)_m=x(x+1)\ldots(x+m-1).
$$

When $\kappa\notin \Z+{1/2}$ the differential 
Galois group is ${\rm SL}_2(\C)$ 
(see e.g. Kolchin [\KOLCHIN]).

See generally e.g. [\BATEMAN, \DLMF, \OLVER].

\bigbreak
\noindent
{\bf Definability}

\medskip

For each fixed  $\kappa>-\tfrac{1}{2}$, consider the   functions
$$
f_1(x)=\sqrt{x^{-1}}e^{-ix^{-1}}H^{(1)}_\kappa(x^{-1}),\quad f_2(x)=\sqrt{x^{-1}}e^{ix^{-1}}H^{(2)}_\kappa(x^{-1}),\quad x>0.
$$

\medbreak
\noindent{\bf 2.1. Lemma.\/}
{\it The functions $f_1$ and $f_2$ extend smoothly to $[0,\infty)$. Moreover, for any $0<\varepsilon<\pi/2$, there is a constant $C=C(\kappa,\varepsilon)>0$ such that
$$
\left|f_j^{(n)}(z)\right|\leq C^{n+1}n!^2,\qquad
\forall n\in\mathbb{Z}_{\geq0},\:\:
0<|z|\leq1,\:\: |\arg z|<\pi-\varepsilon,\:\: j=1,2.
$$
In particular, the real and imaginary parts of $f_1$ and $f_2$ are definable on $[0,1]$ in $\R_{\mathcal{G}}$ and on $[0,\infty)$ in $\R_{\mathcal{G},\exp}$, the o-minimal structures from [\DRIESSPEISSEGGER]}.

\medbreak
\noindent
Holomorphic functions satisfying such derivative bounds on a bounded sector about the positive real axis of opening greater than $\pi$ have Gevrey 1-summable expansions at the origin. In one variable, this is equivalent to the definition in [\DRIESSPEISSEGGER]; see [\BEURLING, \DYNKIN] (cf. [\BALSER, \S4.5 E8] for another equivalent definition).

\medbreak
\noindent
{\bf Proof.\/} We begin with the integral identity, initially for positive real $z$ and with principal branches:
$$
H_{\kappa}^{(1)}(z) = \left(\frac{2}{\pi z}\right)^{1/2} \frac{e^{i(z - \frac{1}{2}\pi\kappa - \frac{1}{4}\pi)}}{\Gamma(\kappa + \frac{1}{2})} \int_0^{\infty} e^{-u} u^{\kappa - \frac{1}{2}} \left(1 + \frac{iu}{2z}\right)^{\kappa - \frac{1}{2}} du
$$
which is stated in [\AWH, \S11.6] for the function $K_\kappa(z)=\tfrac{\pi}{2} i^{\kappa+1}H_\kappa^{(1)}(iz)$.
Thus, with
$$
c_\kappa=\sqrt{\frac{2}{\pi}}\,
\frac{e^{-i(\pi\kappa/2+\pi/4)}}{\Gamma(\kappa+\frac{1}{2})},
$$
we have
$$
f_1(z)=c_\kappa \int_0^{\infty} e^{-u} u^{\kappa - \frac{1}{2}} \left(1 + \frac{iuz}{2}\right)^{\kappa - \frac{1}{2}} du.
$$
If $\kappa=m+\tfrac{1}{2}$ with $m\in\mathbb{Z}_{\geq0}$, this integral is a polynomial in $z$ of degree $m$, so the derivative estimates follow directly. For the derivative estimates below, assume that $\kappa\notin\mathbb{Z}+\tfrac{1}{2}$.

Since
$$
\frac{d^n}{dt^n} (1+t)^{\kappa-\frac{1}{2}} = \frac{\Gamma\left(\kappa + \frac{1}{2}\right)}{\Gamma\left(\kappa - n + \frac{1}{2}\right)} (1+t)^{\kappa - n - \frac{1}{2}},
$$
we get
$$
f^{(n)}_1(z)=c_\kappa\frac{\Gamma\left(\kappa + \frac{1}{2}\right)}{\Gamma\left(\kappa - n + \frac{1}{2}\right)}\frac{i^n}{2^n} \int_0^{\infty} e^{-u} u^{\kappa+n - \frac{1}{2}} \left(1 + \frac{iuz}{2}\right)^{\kappa-n - \frac{1}{2}} du.
$$
Fix $0<\varepsilon<\pi/2$. We use the integration ray
$$
\arg u=
\begin{cases}
-\frac{\pi}{2}+\frac{\varepsilon}{2},
& 0\leq\arg z<\pi-\varepsilon,\\[4pt]
\phantom{-}\frac{\pi}{2}-\frac{\varepsilon}{2},
& -\pi+\varepsilon<\arg z<0.
\end{cases}
$$
These contour rotations are first 
performed for $z$ near the positive real axis, where no branch point or branch cut is crossed. For $u$ on the chosen ray and $z$ in the corresponding sector, we have
$$
\text{Re} u=\sin(\varepsilon/2)|u|,
\qquad
\frac{\varepsilon}{2}\leq\arg(izu/2)\leq\pi-\frac{\varepsilon}{2}.
$$
Thus the rotated integrals give analytic continuations into the indicated sectors, and
$$
|1+izu/2|\geq\sin(\varepsilon/2).
$$
For integers $n>\kappa-\tfrac{1}{2}$, we conclude, after rescaling the ray parameter, that there is a constant $C=C(\kappa,\varepsilon)>0$ such that
\begin{align*}
    \left|f^{(n)}_1(z)\right|&\leq C^{n+1}\frac{\Gamma\left(\kappa + \frac{1}{2}\right)}{\left|\Gamma\left(\kappa - n + \frac{1}{2}\right)\right|}\int_0^\infty e^{-t} t^{\kappa+n - \frac{1}{2}} dt\\
    &=C^{n+1}\frac{\Gamma\left(\kappa + \frac{1}{2}\right)}{\left|\Gamma\left(\kappa - n + \frac{1}{2}\right)\right|}\Gamma(\kappa+n+\tfrac{1}{2}).
\end{align*}
By the reflection formula for the Gamma function and Stirling's approximation, for these $n$ there is a constant $C>0$ such that
$$
\frac{\Gamma(\kappa+n+\tfrac{1}{2})}{\left|\Gamma\left(\kappa - n + \frac{1}{2}\right)\right|}
=\frac{|\cos(\pi\kappa)|}{\pi}
\Gamma(n+\kappa+\tfrac{1}{2})\Gamma(n-\kappa+\tfrac{1}{2})
\leq C^{n+1}n!^2.
$$
For the finitely many integers $0\leq n\leq\kappa-\tfrac{1}{2}$, the same rotated integral gives uniform bounds using $|z|\leq1$ and $|1+izu/2|\leq1+|u|/2$. Enlarging $C$ proves the estimate for all $n$.
The same estimate for the function $f_2$ follows from the identity $H^{(2)}_\kappa(z)=\overline{H^{(1)}_\kappa(\bar{z})}$.
Dominated convergence in the differentiated integrals as $x\to0^+$ also gives the asserted smooth extensions.

The sectorial estimates imply that the real and imaginary parts of $f_1$ and $f_2$ are definable in $\R_{\mathcal{G}}$ near $0$. Their analyticity away from $0$ then gives definability on $[0,1]$. For $x\geq1$, put $t=1/x\in(0,1]$. The convergent expansions of $J_\kappa(t)$ and $Y_\kappa(t)$ express them using functions analytic on a neighborhood of $[0,1]$, real powers of $t$, and, for integer $\kappa$, $\log t$. Since $t^\alpha=\exp(\alpha\log t)$, these expressions are definable in $\R_{\mathrm{an},\exp}$. The factors $t^{1/2}e^{\mp it}$ use only algebraic functions and restricted sine and cosine. Hence the real and imaginary parts of $f_1$ and $f_2$ are definable on $[1,\infty)$ in $\R_{\mathrm{an},\exp}\subseteq\R_{\mathcal{G},\exp}$.
\qed
\bigbreak
\noindent
{\bf Recurrence relation}

\medskip

Bessel functions satisfy 
the following recurrence relations 
(see e.g. [\DLMF, \S10.61])
$$
J_{\kappa+1}(x)={2\kappa \over x} J_\kappa(x)-J_{\kappa-1}(x),
\quad
J_{\kappa+1}(x)=-2J_\kappa'(x)+J_{\kappa-1}(x).
$$
Eliminating $J_{\kappa-1}$ gives a relation between
the functions and derivatives of orders $\kappa,
\kappa+1$.
$$
J_{\kappa+1}(x)={\kappa\over x}J_\kappa(x)-J_\kappa'(x).
$$
Differentiating and using the differential equation 
will give a similar identity for $J_{\kappa+1}'(x)$ 
in terms of $J_\kappa(x), J_\kappa'(x)$, and indeed
(see again [\DLMF, \S10.61]) the same identities
hold for $Y_\kappa (x), H^{(i)}_\kappa(x)$ or
any $\C$-linear combination of them.

\bigbreak
\noindent
{\bf Non-generic order}

\medskip

When $\kappa=m+1/2\in \Z+1/2$, the 
solutions of Bessel's differential 
equation are combinations of
$e^{\pm ix}, \sqrt{x},$ and rational 
functions. Specifically, for a
non-negative integer $m$, the Hankel 
functions are given explicitly by
(see [\DLMF, 10.49(i) via 10.47(ii)])
$$
H^{(1)}_{m+1/2}=(-i)^{m+1}\sqrt{2\over \pi x}
e^{ix}
\sum_{k=0}^mc(m,k) \Big({i\over x}\Big)^k,
$$
$$
H^{(2)}_{m+1/2}=i^{m+1}\sqrt{2\over \pi x}
e^{-ix}
\sum_{k=0}^m c(m,k) \Big({-i\over x}\Big)^k,
$$
where $c(m.k)=(m+k)!/\big((m-k)!k!2^k\big)$ with (as previously)
$$
H^{(1)}_{m+1/2}=J_{m+1/2}+iY_{m+1/2},
\quad
H^{(2)}_{m+1/2}=J_{m+1/2}-iY_{m+1/2}.
$$
For negative non-generic orders one has
$$
H^{(1)}_{-m-1/2}=i(-1)^mH^{(1)}_{m+1/2},
\quad
H^{(2)}_{-m-1/2}=-i(-1)^mH^{(2)}_{m+1/2}.
$$
\bigbreak

\section{Ax-Lindemann for generic 
Bessel equations}

\medskip

This section proves a result on the 
algebraic independence of logarithms 
and solutions of Bessel's differential 
equation with several different generic 
orders on algebraically dependent 
arguments. The initial result, for two 
algebraically dependent arguments, is the
key for the general Bessel-Ax-Schanuel 
statement.

Let $x, y$ be non-constant rational 
functions on a smooth projective complex 
curve $X$ with (as we may assume) 
$\C(X)=\C(x,y)$. We consider the 
differential field $K=\C(x,y)$ with 
derivation $D=d/dx$.

We consider the Bessel d.e. on the $x$ 
variable with generic parameter 
$\kappa$, and on $y$ with generic 
parameter $\lambda$ (and differentiation 
with respect to $y$).

Let $M/K$ be a differential field 
extension, with constant field $\C$, 
that contains subfields $L_x, L_y$, 
where $L_x/K$ is a 
differential extension of $K$ generated 
by a fundamental matrix formed from
two linearly independent solutions 
$z_1(x), z_2(x)$ of Bessel's d.e.
(with order $\kappa$), and $L_y/K$ a
differential extension of $K$ 
generated by a fundamental matrix 
formed by two linearly independent 
solutions $w_1(y), w_2(y)$ of Bessel's 
d.e. (with order $\lambda$).
Thus
$$
{\rm td}_{K}L_x={\rm td}_{K} L_y=3.
$$
Let $L=L_xL_y$. Then $L_x, L_y, L$
are Picard-Vessiot extensions of $K$.

\medbreak
\noindent
{\bf 3.1. Theorem.\/} {\it In the above 
setting, suppose that
$$
{\rm td}_KL < 6.
$$
Then either $y=x$ and 
$\kappa\pm \lambda\in\Z$, or $y=-x$ and
$\kappa\pm \lambda\in\Z$.}

\medbreak
\noindent
{\bf Proof.\/} By the results of Kolchin 
[\KOLCHIN, italicized conclusion at the 
bottom of p 1158], if ${\rm td}_KL < 6$
then there exist $A\in {\rm GL}_2(K)$, $\alpha\in L$,
with $\alpha^2\in K$, such that a relation of the form 
$$
W_y=\alpha A W_x C
$$
holds between any pair of fundamental matrices 
$W_y$, for the Bessel d.e. on $y$, and 
$W_x$ for the Bessel d.e. on $x$,
for some (basis dependent) $C\in {\rm GL}_2(\C)$. 

We take the canonical Hankel function 
solutions on $\PP^1(\C)$, and form
$W_x, W_y$ by composing them with the 
rational functions $x,y$ on the complex 
curve $X$. They are 
related as above, and this relation
persists under analytic continuation
outside of finitely many singular 
points on $X$ (the solutions are in 
general multivalued).

Suppose $P\in X$ is a point in 
${\rm Poleset}(y)$. Then we must 
have also $P\in {\rm Poleset}(x)$,
otherwise the above relation would 
not admit an essential singularity 
in $y$. The converse
also holds, so 
${\rm Poleset}(y)=${\rm Poleset}(x)=
${\rm Poleset}$. 

Let $P\in {\rm Poleset}$. 
Fix a branch of $y$
over $x$ and expand $y$ in a Puiseux 
series in $x$:
$$
y=c_mx^m+c_{m-1/k}x^{m-1/k}+\ldots+
c_1x^{1/k}+{\rm regular\ terms}
$$
where $m\in\Q$ and $k$ is some positive 
integer. The regular terms at $\infty$ 
are those in non-negative powers of 
$x^{-1/k}$. We will denote the sum
of all regular terms at $P$ by 
${\rm Reg}_P(x)$.

We next consider the exponential factors 
in the formal asymptotic expansions of 
the Hankel function solutions at $\infty$
and conclude from Kolchin's relation that, near $P$, we must have
$$
y=\pm x +{\rm Reg}_P(x).
$$

Thus, in particular, $x$ and $y$ have 
the same multiplicity at each point
$P\in {\rm Poleset}$, and hence the same
number of poles (counted with 
multiplicity).



Now if the sign in 
$y=\pm x + {\rm Reg}_Q(x)$ is the 
same at all $Q\in {\rm Poleset\/}$
then $y \mp x$ has no poles and is 
constant. In that case 
$X=\PP^1$ and the zeros of
$x$ and $y$ must coincide, as the 
Bessel functions have monodromy at the origin.
So we have $y=\pm x$.

In this case, the conclusion that $\lambda\pm \kappa\in \Z$
follows from the work of Siegel [\SIEGEL]
(and was reproved by Kolchin [\KOLCHIN]). 
We give another argument.

First, we show that in the relation $W_y=\alpha AW_xC$
we in fact have $\alpha\in \C(x)$. To see this, we observe 
from $\alpha A= W_yC^{-1}W_x^{-1}$
that $\alpha$ is unramified outside $\{0, \infty\}$.
Further, by looking at the formal Hankel 
solutions
at $\infty$, we see that the factors 
$x^{-1/2}$ all cancel, the exponential 
factors match up, and so $\alpha$ is also 
unramified there. Then the divisor
of $\alpha^2$ has even degree at all 
non-zero points, and hence also at $x=0$.
Therefore its square root  $\alpha$
is in $\C(x)$ (whose Jacobian is trivial).

Now suppose $\kappa\notin\Z$.
Consider $W_x$ formed from
$J_{\kappa}(x), J_{-\kappa}(x)$.
Monodromy around $x=0$ is diagonal with 
entries $e(\kappa), e(-\kappa)$. Hence
the monodromy of $W_y$ around 0 has the 
same eigenvalues, so that 
$\{e(\pm\lambda)\}=\{e(\pm\kappa)\}$,
giving $\lambda\pm \kappa \in \Z$.
Otherwise, $\kappa, \lambda$ are both 
integers (so again 
$\lambda\mp \kappa \in \Z$).

Thus to prove the theorem it remains to show that the 
sign in $y=\pm x+{\rm Reg}_Q(x)$ is constant
over $Q\in {\rm Poleset\/}$.


So let us suppose, as we may, that there
is $P\in{\rm Poleset\/}$ at which locally
$y=+x+{\rm Reg}_P(x)$. We take the 
principal Hankel function bases
to form $W_y, W_x$. We denote by
$\tilde{W}_y, \tilde{W}_x$ the 
corresponding formal solutions given 
by the asymptotic expansions
in the right half-plane.

The Kolchin relation
$$
\tilde{W}_y=\alpha A \tilde{W}_x C
$$
is preserved by the summation process.
The exponential torus
$$
T_\mu=
\begin{pmatrix}\mu & 0 \cr
0 & \mu^{-1}\end{pmatrix}, \quad 
\mu\in\C^\times
$$
acts by
$$
\tilde{W}_y \mapsto \tilde{W}_yT_\mu,
\quad
\tilde{W}_x \mapsto \tilde{W}_xT_\mu,
$$
and this implies that $C$ commutes with
$T_\mu$ and is therefore diagonal.

Now suppose $Q\in{\rm Poleset\/}$ 
(with $Q=P$ allowed). We analytically 
continue the relation
$$
W_y=\alpha A W_x C
$$
to $Q$, where (by doing some extra loop
in $x$ near $Q$) we may assume that $W_x$
has returned to the principal branch,
while $y$, in the limit of 
$x\rightarrow\infty$, makes $m$ 
``half-loops'' by $\pi i$ where $m\in\Z$ as we have 
$y=\pm x+{\rm Reg}_Q(x)$ near $Q$.
The plus sign corresponds to even $m$,
the minus sign to odd. 

We work on the Riemann surface of the logarithm
and write $y=y^*e^{m\pi i}$ where
$y^*$ is in the principal branch. So $y^*= x +{\rm Reg}_Q(x)$.

We now apply the ``analytic continuation''
(with reflection) formula from
[\DLMF, \S10.11]:
$$
W_{ze^{m\pi i}}=N_mW_{z}S_{\nu, m}
$$
for the Hankel functions of order $\nu$
where
$$
N_m=\begin{pmatrix}1 & 0 \cr
0 & e^{m\pi i}\end{pmatrix}, 
$$
$$
S_{\nu, m}={1\over \sin(\nu\pi)}
\begin{pmatrix}-\sin((m-1)\nu\pi)
& 
e^{\nu\pi i} \sin(m\nu\pi)
\cr
-e^{-\nu\pi i} \sin(m\nu\pi)
& 
\sin((m+1)\nu\pi)\end{pmatrix},
$$
and where the limiting values are taken
for integral $\nu$, namely, for
$\nu\in\Z$,
$$
S_{\nu, m}=
\begin{pmatrix}
-(-1)^{m\nu}(m-1) 
& 
(-1)^{(m-1)\nu}e^{\nu\pi i} m
\cr
-(-1)^{(m-1)\nu} e^{-\nu\pi i} m
& 
(-1)^{m\nu}(m+1)\end{pmatrix}. 
$$

We have
$$
W_{y^*e^{m\pi i}}=
N_mW_{y^*}S_{\lambda, m}=
\alpha A W_x C,
$$
so that
$$
N_m\tilde{W}_{y^*}S_{\lambda, m}=
\alpha A \tilde{W}_x C,
$$
and, under the exponential torus,
$$
N_m\tilde{W}_{y^*}T_\mu S_{\lambda, m}=
\alpha A \tilde{W}_x T_\mu C.
$$
Comparing the last two relations, 
and since $C^{-1}T_\mu C=T_\mu$, 
we find
$$
S_{\lambda, m}^{-1}T_\mu S_{\lambda, m}
= C^{-1}T_\mu C=T_\mu.
$$
This entails that $S_{\lambda, m}$
is a diagonal matrix.

If $\lambda$ is irrational then this is 
only possible if $m=0$, which immediately
gives that $y=+x +{\rm Reg}_Q(x)$.
If $\lambda\in \Z$ we are again reduced to
$m=0$.

So we are left with $\lambda=a/b\in\Q$
for relatively prime integers $a,b$
with $b>2$ as $\lambda\in\Z+1/2$ is
non-generic and excluded.
The condition that $S_{\lambda, m}$ is
diagonal now entails that 
$\sin(m\lambda\pi)=0$ and hence $b|m$. 

Let us consider the
(possibly multivalued) function
$$
\phi=\Big({y\over x}\Big)^{2/b}
$$
on $X$, normalized on some branch by 
$\phi(P)=1$. Note that $\phi$ is only ramified at the 
zeros and poles of $y/x$, and in particular not in ${\rm Poleset}$.

As any path from $P$ to $Q$ may be locally put in the above
form (by making suitable loops around $Q$),
the condition $b|m$ implies 
that $\phi$ is single-valued at $Q$, and 
indeed in a neighbourhood of $Q$, as 
multiple values would differ by roots 
of unity multiples and so  differ too 
much in their argument near $Q$. 
So $\phi$ has only one branch near $Q$
and hence only one branch
(indeed one point $Q\in{\rm Poleset}$ with this 
property would suffice for
this). Thus,
$\phi$ is single-valued on $X$. 

Since $b|m$, we see that $\phi$ takes 
the value $1$ at all 
$Q\in {\rm Poleset\/}$. If, at $Q$,
$x, y$ have (the same) multiplicity $k$
then $\phi$ takes the value $1$ at $Q$ with
multiplicity $k$. Hence $\phi$ takes 
the value $1$ with total multiplicity
at least $N$, where $N$ is the degree of $x$. This equals the degree of $y$, as $x,y$ have the same poles with the same multiplicities. On the other hand
$\phi$ has its zeros only at points
where $y$ is zero. 
The total zero multiplicity of $\phi$
is then at most $2N/b$, which is 
strictly smaller.

Hence $\phi$ is constant, and we conclude
that $y = x$. \qed

\medbreak

The transcendence indeed drops in the 
exceptional cases.

Next we consider non-constant rational 
functions
$$
e_1,\ldots, e_m, x_1,\ldots, x_n
$$
on a smooth projective complex curve $X$ (with $\C(X)=\C(e_i, x_j)$).
We have (multivalued) functions
$t_j=\log e_j$ and fundamental matrices $W_i$ (whose entries are multi-valued functions)
for Bessel's differential equation of 
generic order $\kappa_i$ in the 
argument $x_i$.

\medbreak
\noindent
{\bf 3.2. Theorem.\/}
{\it In the above setting, suppose that
$$
{\rm td}_{\C(x_1)}\Big(\overline{e}, 
\overline{t},
\overline{x}, \overline{W}\Big) < m+3n.
$$
Then either (A) or (B) holds where:

(A) The $t_j$ are $\Q$-linearly 
dependent modulo $\C$.

(B) There exists distinct indices $j,k$ 
of the generic Bessel equations such that 
either

(i) $x_k=x_j$ (and $\kappa_j\pm \kappa_k 
\in \Z$), or

(ii) $x_k=-x_j$ (and $\kappa_j\pm \kappa_k \in \Z$).\/}

\medbreak
\noindent
{\bf Proof.\/}
By Kolchin's results [\KOLCHIN: Theorem p 1152, Ostrowski's theorem p 1155, and italicized text bottom p 1158] 
applied to the Picard-Vessiot extension of $\C(X)$ we get either:

(a) a $\C$ linear dependence modulo 
$\C$  of the $t_j$, but then by 
Ax-Schanuel [\AX] this must be a 
$\Q$-linear dependence mod $\C$, 
giving (A), or

(b) a pair $i,j$ with $i\ne j$ and 
a relation as in Theorem 3.1. 
This leads to (B).\ \qed


\section{Weakly special subvarieties 
and Ax-Schanuel}

\medskip

We start with the main result of 
[\BCFNDAS]. The setting is a principal 
$G$-bundle $P\rightarrow Y$ where $Y$ 
is a complex algebraic variety and $G$ 
is a complex algebraic group.

\medbreak

\noindent
{\bf 4.1. Theorem.} ([\BCFNDAS, Theorem A])
{\it Let $\nabla$ be a $G$-principal flat connection on $P\rightarrow Y$
with sparse Galois group ${\rm Gal}(\nabla)=G$.
Let $V$ be an irreducible algebraic subvariety of $P$, $\mathcal{L}$ a horizontal leaf,
and $A\subset V\cap\mathcal{L}$ an irreducible analytic subvariety of $P$.

If $\dim V < \dim G + \dim A$ then the 
projection of $A$ in $Y$ is contained 
in a $\nabla$-special subvariety. \qed}

\medbreak

We refer to [\BCFNDAS, Lemma 2.3] for 
the definition of the Galois group 
${\rm Gal}(\nabla)$ but note that 
$\nabla$-special subvarieties, 
which for us are the (strict) 
maximal weakly special subvarieties, 
are the irreducible maximal subvarieties 
$X\subset Y$ such that
${\rm Gal}(\nabla\vert_X)$ is a strict 
subgroup of ${\rm Gal}(\nabla)$.
See [\BCFNDAS, Definition 3.4 (and 
informally at p1)] for the definition 
of a ``sparse'' algebraic group. We note 
that powers of $\GG_{\rm m}$ are not 
sparse, but products of $\GG_{\rm a}$ 
and ${\rm SL}_2$ are sparse. Hence we 
work with the logarithm rather than the 
exponential.

We consider solutions of a system of logarithm and Bessel differential equations.
We take $m$ independent copies of the 
logarithm d.e. $dt_j/de_j=1/e_j$
and $n$ copies of the Bessel d.e. with 
generic orders $\kappa_1,\ldots, \kappa_n$
on variables $x_1, \ldots, x_n$.
So we have $G=\GG_{\rm a}(\C)^m \times 
{\rm SL}_2(\C)^n$ and a principal
$G$ bundle $P\rightarrow Y$, with a 
$G$-principal flat connection $\nabla$, 
where $
Y= (\C^\times)^m \times (\C^\times)^n,
$
and taking $\overline{t}=(t_1,\ldots, 
t_m)$ and $\overline{h}=(h_1,\ldots, 
h_n)$, with
$$
P=\{(\overline{t}, \overline{h},
\overline{e}, \overline{x})
\in \big(\GG_{\rm a}(\C)\big)^m \times 
\big( {\rm GL}_2(\C)\big)^n
\times(\C^\times)^m \times (\C^\times)^n:
$$
$$
\det(h_j)={2\over \pi x_j}, 
j=1,\ldots, n\}.
$$
Locally on $Y$, the field generted by
fundamental matrices of solutions has
transcendence degree $m+3n$
over $\C(Y)$ as the variables are
independent.

To get an explicit Ax-Schanuel statement 
we must explicate the maximal weakly 
special subvarieties of $Y$. In this 
Picard-Vessiot setting, they are the 
maximal subvarieties $W\subset Y$
where the differential Galois group 
is a strict subgroup of $G$
(see [\BCFNDAS, Examples 2.1 and 2.5]). 

\medbreak
\noindent
{\bf 4.2. Proposition.\/} {\it In this setting, the maximal weakly special subvarieties
$W\subset Y$ are the hypersurfaces of the following types:

\smallskip

(i) $W$ is defined by an equation 
$e_1^{k_1}\cdot\ldots\cdot e_m^{k_m}=c$,
where $k_i\in \Z$ are not all zero, 
and $c\in \C^\times$;

\smallskip

(ii) $W$ is defined by $x_i=c$ 
for some $i$ and some 
$c\in \C^\times$;

\smallskip

(iii) $W$ is defined by
$x_j=x_k$ for some $k,j\ne k$ with 
$\kappa_k\pm \kappa_j\in \Z$; 

\smallskip

(iv) $W$ is defined by 
$x_j=-x_k$ for some $k,j\ne k$ with 
$\kappa_k\pm \kappa_j\in \Z$.}

\medbreak
\noindent
{\bf Proof.}
If a subvariety $W\subset Y$ is 
not contained in a hypersurface 
of such type then it contains a 
curve $C$  not contained in any
such hypersurface, and on such 
a curve the transcendence degree 
of the extension does not drop. 

To see this, suppose $W\subset Y$ 
is a positive dimensional variety 
on which none of the following countably infinite list of conditions
holds:

(i)${}_{\overline{k}}$: 
$e_1^{k_1}\cdot\ldots\cdot 
e_m^{k_m}=c$,
where $k_i\in \Z$ are not all zero, 
and $c\in \C^\times$;

\smallskip

(ii)${}_i$: $x_i=c$ 
for some 
$c\in \C^\times$;

(iii)${}_{jk}$: $x_j=x_k$ 
where $j\ne k$ are such that 
$\kappa_k\pm \kappa_j\in\Z$; 

(iv)${}_{jk}$: $x_j=-x_k$ 
where $j\ne k$ are such that 
$\kappa_k\pm \kappa_j\in\Z$.

\smallskip

The claim is that $W$ contains an
algebraic curve on which none of 
these conditions hold.

To prove the claim, let $P\in W$ 
be a smooth point with 
all coordinates non-zero. We consider
codimension 1 subspaces of the tangent space to $W$
at $P$. For each of our conditions, 
only finitely many such subspaces will
give the corresponding condition 
holding locally. 

So we can find a direction avoiding 
all these conditions, and find 
an algebraic curve  $C\subset W$ 
through $P$ with that local direction.
\ \qed

%
%

\medbreak
Now Theorem 4.1 gives a complex form of Ax-Schanuel for Bessel and logarithm functions.
It is then a standard matter to deduce a differential field version of Ax-Schanuel from
the complex version.
We can include Bessel d.e. of non-generic order since
$\C(e^{iy}, iy)$ is a finite extension of $\C(y, W(y))$
when $W(y)$ is a fundamental matrix of a Bessel d.e. of 
non-generic order.

\medbreak
\noindent
{\bf 4.3. Proof of Theorem 1.4.\/}  
First, we get rid of the non-generic Bessel d.e. by including further logarithm d.e,
and renumber so that we work with $e_j$, their logarithms (i.e. logarithm d.e. solutions) $t_j$,
$x_j$ and fundamental matrices for generic Bessel d.e. of orders $\kappa_j$.
We have assumed all $e_j$ and $x_j$ are non-constant. Taking a Seidenberg embedding,
we have local complex solutions, hence some irreducible one dimensional
analytic subvariety $A$ of a horizontal leaf of
the bundle $P$. Let $V$ be its Zariski closure. Possibly enlarging $A$ we can assume $A$ is
a component of the intersection of $V$ with a horizontal leaf of the connection. If
$$
{\rm td}_{\C}\C\Big(\overline{t}, \overline{e}, \overline{x}, 
\overline{W}\Big) < \ell+3n+1
$$
then (as ${\rm rank} A=1$ and $\dim G =\ell + 3n$) we have
$$
\dim V< \dim G + \dim A
$$
and we conclude that the locus of $e_i, x_j$ is in a proper weakly special subvariety,
and then the constants defining it are in $\C$. \ \qed

\medbreak

One could also give a version with several commuting derivations.


\section{Two Lemmas}

In this section, $T=\{(t_1(s),\ldots, 
t_n(s)): s\in \Delta\}$ 
is a local complex algebraic 
curve, meaning that we consider the 
functions $t_i(s)$ to be defined on 
some open ball $\Delta\subset \C$ 
(say centered at $s=0$). We assume
the image $T$ is a relatively open
subset of a complex algebraic curve.
We may assume without loss of 
generality that all the $t_i(s)$ 
are holomorphic for $s\in \Delta$.

The proofs involve only local arguments
(in particular, appeals to Ax-Schanuel),
except that in the proof of Lemma 5.2 
we eventually reduce to global 
considerations concerning a certain 
rational function.

\medskip
\noindent{\bf 5.1. Generic Lemma.\/}
{\it Let $\kappa\in \C$ be generic.
Let $V\subset (\C^\times)^n$ be an 
algebraic subvariety.
Suppose $T=(t_1(s),\ldots, t_n(s))$,
where $s\in \Delta$,
is a local complex algebraic 
curve in $\C^n$ and that 
$\overline{x}(s)\in V$ where
$\overline{x}(s)=\big(x_1(s),\ldots, x_n(s)\big)$ 
is analytic for $s\in\Delta$ and satisfies
$$
e^{-2\pi i t_j}H^{(1)}_\kappa(x_j)+e^{2\pi i t_j}
H^{(2)}_\kappa(x_j)=0, \quad j=1,\ldots, n.
$$
Then there is a weakly special subvariety $W\subset V$
containing $\{\overline{x}(s)\}$.\/
}

\medbreak
\noindent
{\bf Proof.\/} We suppose 
$V\subset (\C^\times)^n$
is strict or the conclusion holds with 
$W=V=(\C^\times)^n$.
Let us observe that the
system of functions
$$
2\pi i t_j, \quad e^{2\pi i t_j},\quad  
x_j, \quad W_j,
$$
where $W_j$ is a fundamental matrix for 
Bessel's differential equation with parameter $\kappa$,
has transcendence degree at most
$$
1+n+\dim V + 2n \le n+ 3n
$$
(as $\dim V<n$)
and thus we are in the exceptional case of
Ax-Schanuel. Let us suppose without loss of generality
that $x_1,\ldots, x_k$ are constant, the other
$n'=n-k$
coordinates $x_j$ (and corresponding $t_j$)
non-constant. And suppose $m''$ of the
non-constant $t_j$ are $\Q$-linearly independent
modulo $\C$. Suppose further that $n''$ of the
non-constant
$x$ coordinates are maximal with $x_j\ne \pm x_i$.
Let $W$ be the weakly special subvariety defined
by the constant coordinates and equalities and
anti-equalities among the $x_i$, and let $V''\subset W\cap V$ be the irreducible component containing $\overline{x}(s)$.
Then $\dim W= n''$.
Considering the non-constant $\pm$-distinct
functions $x_j$ and corresponding $t_j$
we have transcendence degree at most
$$
1+m''+\dim V'' + 2n'' 
$$
If $\dim V''< \dim W=n''$ we are again exceptional
for Ax-Schanuel, which implies a further weakly
special relation and
contradicts our definition
of $W$. So $W\subset V$ and 
$\{\overline{x}(s)\}\subset W$.\ \qed

\medbreak

Now we establish the corresponding 
result for a non-generic order $\kappa$, 
finding that the specific form of the 
Hankel functions entails exactly the 
same conclusion as the generic case.

\medbreak
\noindent{\bf 5.2. Non-Generic Lemma.\/}
{\it Let $\kappa\ne \pm 1/2$ be non-generic.
Let $V\subset (\C^\times)^n$ be an 
algebraic subvariety.
Suppose $T=(t_1(s),\ldots, t_n(s))$,
where $s\in \Delta$,
is a local complex algebraic 
curve in $\C^n$ and that 
$\overline{x}(s)\in V$ where
$\overline{x}(s)=\big(x_1(s),\ldots, x_n(s)\big)$ 
is analytic for $s\in\Delta$ and satisfies
$$
e^{-2\pi i t_j}H^{(1)}_\kappa(x_j)+e^{2\pi i t_j}
H^{(2)}_\kappa(x_j)=0, \quad j=1,\ldots, n.
$$
Then there is a weakly special subvariety $W\subset V$
containing $\{\overline{x}(s)\}$.\/
}

\medbreak
\noindent
{\bf Proof.} As in 5.1 we may suppose $V$
is strict.
In the first part of the
proof, we work only locally.
We may assume that 
$\kappa$ is positive. We may also 
assume that $V$ is a strict
subvariety or the conclusion is 
immediate. The given relation leads to
(locally)
$$
e^{4\pi i t_j}=R(x_j)e^{2i x_j}
$$
for $j=1,\ldots, n$ where $R$ is a rational function
$R=R_\kappa$ where $\kappa=m+1/2$ for 
some positive integer $m$, given by
$$
{H^{(1)}_\kappa(x)\over H^{(2)}_\kappa(x)}=
-R(x) e^{2ix},
$$
specifically (see \S2)
$$
R=
(-1)^{m}
{\sum_{k=0}^m c(m,k)(i/x)^{k} \over 
\sum_{k=0}^m c(m,k)(-i/x)^{k}
}.
$$
We observe that $R$ is non-constant (as 
$\kappa\ne \pm 1/2$) and
has modulus 1 on $\R^\times$
(as numerator and denominator are complex conjugates),
and moreover $R(x)R(-x)=1$.

Now some $x_j$ may be constant, and are so if and only if the
corresponding $t_j$ is constant. Renumbering
we ignore these and again conclude if $V$ is not
a strict subvariety on the non-constant coordinates.

So we may and now do assume that all 
coordinates are non-constant, and $V\subset (\C^\times)^n$ is strict.

Let us suppose that the $\Q$-linear mod $\C$ dimension
of $\{4\pi it_j\}$ is $d$, that the $\Q$-linear mod $\C$ 
dimension of $\{4\pi i t_j, 2i x_h\}$ is $d+\delta$,
and that the Zariski closure
of the locus of $x_j$ in $(\C^\times)^n$
has dimension $k$ (with $k\le n-1$ as $V$ is a strict subvariety).

Write $y_j=R(x_j)$ noting that, locally,
$$
\log y_j=4\pi i t_j -2i x_j.
$$
Note also that, since $R$ is non-constant,
${\rm td}_{\C}\C(\overline{y})=
{\rm td}_{\C}\C(\overline{x})=k$.
Let further $r$ be the multiplicative rank of $\{y_j\}$ modulo constants.

By Ax-Schanuel, taking logarithms 
locally, we have
$$
1+r \le {\rm td}_{\C}\C\Big(
\log \overline{y}, \overline{y}\Big).
$$
The $t_j$ are pairwise algebraic over each other. So an upper bound for
this transcendence degree is $1+k$.
If however the $t_j$ are algebraic over
$\C(\overline{x})$ then transcendence
degree is $k$ and $r<k$ is a 
contradiction.
So we conclude that the $t_j$ are not
algebraic over $\C(\overline{x})$,
and $r=k$.

Now we observe that
$$
{\rm td}_{\C}\C\Big(4\pi i\overline{t}, 
\exp\big(4\pi i \overline{t}\big), 2i\overline{x},
\exp\big(2i\overline{x}\big)\Big)
\le \min(1+\delta+d+0, 1+k+d+0)
$$
(as each $\exp(2i x_j)$ is  algebraic over
$x_j, \exp(4\pi i t_j)$), 
the minimum permitted by Ax-Schanuel.
Therefore
$$
1+\delta \le {\rm td}_{\C}\C\Big(
4\pi i\overline{t}, 2i\overline{x}\Big)
\le 1+\delta+d-d
$$
and we have equality. But this transcendence degree is $1+k$, and we
conclude that $k=\delta$.

%


Now suppose (after renumbering) that
$2ix_1,\ldots, 2ix_k$ are $\Q$-linearly independent over the $4\pi i t_j$ modulo
constants. The remaining $2i x_h$
satisfy $\Q$ linear relations mod constants over $2ix_1,\ldots, 2ix_k$
and the $4\pi i t_j$. But as the $t_j$
are not algebraic over $\overline{x}$,
we see that these relations must
be just $\Q$-linear mod constant
over $2ix_1,\ldots, 2ix_k$.


Therefore, all the algebraic relations (and there 
is at least one) defining the Zariski closure
of the locus $\{\overline{x}\}$ in $(\C^\times)^n$
are $\Q$-linear-mod-$\C$ relations.
We have that $x_1,\ldots, x_k$
are independent while
$$
x_h=\sum q_{hj} x_j+c_h, \quad h=k+1,\ldots, n
$$
where $q_{hj}\in \Q$ are not all zero and $c_h\in \C$.


Now also $y_1, \ldots, y_k$
are algebraically independent over $\C$, 
and defined by multiplicative relations
modulo $\C$, so we conclude
that there are multiplicative relations
$$
y_h^{\alpha_{hh}}=b_h\prod_{j=1}^k y_j^{\alpha_{hj}}, \quad
h=k+1,\ldots, n,
$$
where $\alpha_{hh}\ne 0$ and
$\alpha_{hj}\in\Z, j=1,\ldots, k$ 
are not all zero, and $b_h\in \C$. 

Now we can work globally with the
two algebraic varieties and the rational map.
Combining, we have (for $h=k+1,\ldots, n$)
$$
R\Big(\sum_{j=1}^k q_{hj}x_j+c_h\Big)^{\alpha_{hh}}=
b_h\prod_{j=1}^k R(x_j)^{\alpha_{hj}}.
$$

As $R$ is non-constant with its 
poles conjugate to its zeros, and 
vice versa, it has at least one pole. 
Then since the variables $x_1,\ldots, x_k$ are
independent,
for each $h$ there can be at most one non-zero $\alpha_{h\ell}$, 
and by degree considerations,
and it can only be $\alpha_{h\ell}=\pm \alpha_{hh}$.  Then
$$
\sum_{j=1}^k q_{hj}x_j+c_h = q_{h\ell} x_\ell+c.
$$ 

So we are reduced to considering (dropping subscripts,
and according as $\alpha=\pm 1$) the possibilities
$$
R(qx+c) =bR(x), \quad R(qx+c)R(x)=b
$$

We consider the first ($\alpha=1$) case.

Suppose $q=1$. This is impossible if $c\ne 0$,
while if $c=0$ we have $b=1$ and equality holds.

Suppose $q\ne 1$. Consider the map $x\mapsto qx+c$.
If we write $\xi = x-\delta$ then the transformation
on $\xi$ is $\xi \mapsto q\xi +(q\delta+c-\delta)$.
So since $q\ne 1$ we can choose $\delta$ such that
the map
on $\xi$ is $\xi\mapsto q\xi$. If this map has
infinite order then such a relation 
$R^*(q\xi)=bR^*(\xi)$, where $R^*(\xi)=R(\xi+\delta)$ 
is impossible as $R^*$ has a zero (and a pole).

This leaves us to consider $q=-1$. As $x\mapsto -x+c$
is an involution we must have $b^2=1$ so $b=\pm 1$.
Also we know that $R(x)R(-x)=1$. So
$$
R(-x+c)=bR(x)=b/R(-x) 
\quad{\rm therefore} \quad  R(-x)R(-x+c)=b
$$
which we may rewrite as
$$
R(x) R(x+c)=b=\pm 1.
$$
Now if $x_0$ is a zero of $R$ then $x_0+c$ is a pole,
while if $x_0+c$ (in first factor) is a pole then
$x_0+2c$ is a zero. This gives a contradiction
if $c\ne 0$. But $c=0$ is also impossible as $R$
is non-constant.

So for $\alpha=1$ the only possibility is 
the trivial $q=1, c=0, b=1$, i.e. $x_h=x_j$ for
some $j=1,\ldots, k$.

Now we consider $\alpha=-1$, i.e. $R(qx+c)R(x)=b$.

If $q=1$ this is impossible as above. 
If $q\ne \pm 1$ then the map $x\mapsto qx+c$ has
infinite order (recall $q\in \Q$) and leads to a
contradiction.

So we consider $q=-1$ and we have the relation
$R(x)R(-x)=1$ so we need to consider if a second
such relation (with $c\ne 0$) is possible. It isn't,
as the orbits under the two maps are infinite.

So for $\alpha=-1$ the only possibility is 
$R(x)R(-x)=1$, i.e. $x_h=-x_j$ for some 
$j=1,\ldots, k$.

And such a relation must occur unless $V$ contains
the weakly special given by the constant $x_i$.

Repeating on the new weakly special we eventually
find a positive-dimensional weakly special subvariety
$W\subset V$ containing the locus $\{\overline{x}(s)\}$.
\ \qed


\section{Bessel-Manin-Mumford}

\medbreak
\noindent
{\bf 6.1. Proof of Theorem 1.2.\/} 
Let $\kappa\in [0, \infty)\backslash\{ 
1/2\}$ 
and $V\subset (\C^\times)^n$.

We can deal with the finitely many 
non-real zeros by a simple reduction
after which it suffices to consider 
non-zero real zeros only
(which, as has been noted, 
are symmetrically located 
about $x=0$).

We prove the theorem by induction on 
the ambient dimension $n$. 
The conclusion is clear
for $n=1$. Now we consider 
$V\subset (\C^\times)^n$ and suppose 
the theorem holds for all subvarieties
$U\subset (\C^\times)^m, m<n$.

The weakly special subvarieties in 
$(\C^\times)^n$ come in finitely 
many families, indexed by the subset
of variables that is constant (which may be empty) and the equalities (or anti-equalities)
among other coordinates. There are only finitely many weakly special subvarieties that have 
no constant coordinates. If $V$ is strict, we need consider only such ``strongly special'' weakly specials
with some pair of equal coordinates, and the conclusion for these follows by induction.
For any other family $T$ of weakly special subvarieties of positive 
dimension, the set of constant coordinates has size $m(T)<n$ and the weakly special subvarieties
of this type that are contained in 
$V$ constitute a subvariety $V_T\subset 
(\C^\times)^{m(T)}$. The special 
subvarieties
of type $T$ contained in $V$ correspond 
to special points in $V_T$
By induction the conclusion holds for 
each $V_T$.

The union of all positive-dimensional
weakly special subvarieties of $V$ is
an algebraic subvariety $V^{\rm pwsp}$, and it suffices to show
that the estimate in the conclusion of Theorem 1.2 holds for $V-V^{\rm pwsp}$,
in which all special subvarieties are special points.

We have the definable set 
$$
Z=Z_{n, \kappa}(V)=\Big\{(t_1,\ldots, 
t_n)\in \R^n: \exists x_1, \ldots, x_n
\in (\R^\times)^n: 
$$
$$
(x_1,\ldots, x_n) \in V-V^{\rm pwsp}, \quad\quad
0\le |x_j|-2\pi t_j < 2\pi, \quad j=1,\ldots, n, \quad {\rm and}
$$
$$
e^{-2\pi i t_j}H^{(1)}_\kappa(|x_j|)+
e^{2\pi i t_j}H^{(2)}_\kappa(|x_j|)=0, 
\quad j=1,\ldots, n\Big\}.
$$

Given $\overline{t}\in Z$, by the definability, 
there is a uniform finite
bound on the number of $\overline{x}$
satisfying the condition (which cannot
hold on an interval of positive length).

Any locally defined positive-dimensional 
semi-algebraic curve in $Z$,
say $\big(t_1(s),\ldots, t_n(s)\big)$
for $s\in I$ some real interval,
may be locally complexified.
Further, the corresponding
functions $x_j(s), s\in I$ may be 
assumed to be real analytic
(perhaps on a smaller interval), by 
analytic cell decomposition in 
$\R_{\rm an}$, as the implicitly defined 
function $x(t)$ is locally analytic,
and hence may also be complexified
(perhaps on a smaller disk $\Delta$).
Moreover, if $t_i(s)$ is non-constant 
then $x_i(s)$ is non-constant 
(and vice versa). Locally, we have
$|x_j|=\sigma_j x_j$ with a locally constant
sign $\sigma_j\in\{\pm 1\}$. 
So, by the Lemmas in \S5,
$(|x_1|,\ldots, |x_n|)$ lies
in a positive-dimensional weakly special
subvariety $W$ of $V^{\overline{\sigma}}=
\{(x_1,\ldots, x_n): (\sigma_1x_1,\ldots, 
\sigma_nx_n)\in V\}$. Hence $\overline{x}$
lies in the positive-dimensional
weakly special subvariety 
$W^{\overline{\sigma}}\subset V$.

Now we apply block counting to $Z$ and 
conclude  that the integer points in $V$ 
are contained in at most 
$c(Z, \epsilon)T^\epsilon$ blocks. 
For blocks of positive dimension,
the corresponding points in $V$ would 
lie in some 
positive-dimensional weakly special $W\subset V$, 
which we have removed. So all the blocks are 
zero-dimensional.


There are at most 
$c(Z, \epsilon)T^\epsilon$ 
zero-dimensional 
blocks, each corresponding to 
a uniformly finite number of isolated
special points in $V$. The proof is therefore complete 
by induction. 

When $\kappa$ is non-generic (but 
$\kappa \ne 1/2$), all the counting 
estimates are poly-logarithmic by [\BNZ].  \qed

\medskip
\noindent
{\bf 6.2. Remark.\/} We note that all constants are 
uniform for $V$ in a  definable family. Moreover,
as Bessel zeros for real order are real (with finitely 
many exceptions) we can use the uniform definability
of real power functions to achieve uniformity in fewnomial
families. This is simpler than the corresponding
deduction in [\BKPLMM] where a fewnomial differential 
elimination result is required as well. 


\section{Bessel distance sets and Fuglede problem}

\medskip

Up to dilation, the Fuglede problem concerns sets $\Lambda\subset \R^m$ whose
non-zero distances are zeros of $J_{m/2}$. We will call such a set a Bessel distance set,
even more generally when the non-zero distances of points in $\Lambda\subset \R^m$
are zeros of $J_\kappa$ for some real number $\kappa$ that is not necessarily equal to
$m/2$. We exclude $\kappa= \pm 1/2$ for which the distances are (up to scaling) integers.
We will prove the following, which clearly implies Theorem 1.5.

\medskip
\noindent
{\bf 7.1. Theorem.\/} {\it
Let $\kappa\in \R_+, \kappa\ne 1/2, m\ge 1, \epsilon>0$ and 
$\Lambda \subset \R^m$ a set of points whose
non-zero distances are all zeros of $J_\kappa$. Then there is a constant $C=C(\kappa, m, \epsilon)$
such that, for $R\ge 1$,
$$
\hash \Lambda \cap [-R, R]^m \le C R^\epsilon.
$$
If $\kappa=r+1/2$ for an integer $r\ge 1$ this estimate may be improved to
$$
\hash \Lambda \cap [-R, R]^m \le C\big(\log R\big)^c.
$$
for suitable constants $C=C(m, \kappa), c=c(m, \kappa)$ and $R\ge e$.
}

\medbreak

Given positive integers $k,m$ one can consider the {\it distance variety\/} $V_{k, {\R^m}}$
defined to be the image in $\R^{k(k-1)/2}$ of the map
$$
\big(\R^m\big)^k\rightarrow \R^{k(k-1)/2},
$$
$$
\big(P_1,\ldots, P_k\big) \mapsto \big(d(P_i, P_j): 1\le i < j\le k\big)
$$
on points $P_1, \ldots, P_k\in \R^m$.
Here we do fix some ordering of the pairs $(i,j)$, so $V$ is not just the variety
of distances but retains a little more information. The image is invariant under 
permuting the points. It is well known that if $k\ge m+2$ then the image is contained
in a strict subvariety; in particular, if $k=m+2$, the {\it Cayley-Menger determinant\/}
$$
{\rm CM}_{m}(\overline{d})=
\begin{pmatrix}
0 & d_{12}^2 & d_{13}^2 & \ldots & d_{1k}^2 & 1 \cr
d_{12}^2 & 0 & d_{23}^2 & \ldots & d_{2k}^2 & 1 \cr
d_{13}^2 & d_{23}^2 & 0 & \vdots & d_{3k}^2 & 1\cr
\vdots & \vdots & \vdots & \ddots & \vdots & \vdots \cr
d_{1k}^2 & d_{2k}^2 & d_{3k}^2 & \ldots & 0 & 1\cr
1 & 1 & 1 & 1 & 1 & 0
\end{pmatrix}
$$
vanishes on the image. (The image is also governed by suitable triangle inequalities.)

We will see that $V_{4, \R^2}$ contains 
no positive-dimensional weakly Bessel 
special subvarieties, other than 
degenerate ones in which some points
coincide. Hence Bessel-Manin-Mumford 
implies something notionally stronger 
than Theorem 7.1, that the number of 
4-tuples of points in $\Lambda$ whose 
distances are in $[R, R]^6$ is 
$R^\epsilon$.

But $V_{5, \R^3}$ does contain positive 
dimensional non-degenerate weakly
special subvarieties. For example, if 
$4$ points, considered fixed, lie on a 
circle in a 2-plane then a point on the normal to the 2-plane at the center of 
the circle will be equidistant to
all four points and its distance can 
vary. The distance has some minimum (the 
radius of the  circle), and so formally 
$V_{4, \R^2}$ does not contain all the 
real points of the corresponding
weakly Bessel special subvariety. 
Now Fuglede's conjecture would 
exclude the existence of such a four 
tuple with Bessel distances. That 
however seems to require proving their
algebraic independence and seems out 
of reach.

We will then work with a variety formed 
by fixing some initial points 
$P_2,\ldots, P_{m+2}$ in suitably general 
position and considering the image of
$$
\R^m \rightarrow \R^{m+1}, \quad P_1 \mapsto \big(d(P_1, P_2),\ldots, d(P_1, P_{m+2})\big).
$$

Further, to avoid awkward issues about weakly specials being contained ``in full dimension''
yet not fully contained in varieties under consideration (as above), it is convenient to
work with complex algebraic varieties.

Thus, given $P_2,\ldots, P_{m+2}\in \C^m$ we consider the variety
$V_{\overline{P}}\subset \C^{m+1}$ defined by the equation
$$
{\rm CM}_{m}(\overline{d})=0.
$$

In what follows, by points in $\R^M$ being in {\it affine general position\/} we mean 
that no $m+1$ of them lie on any $m-1$-dimensional affine linear subspace. 

\medbreak
\noindent
{\bf 7.2. Proposition.\/} 
{\it If the $m+1$ points $P_2,\ldots, P_{m+2}\in \R^m$ are in affine general position 
then any weakly  special subvariety of positive dimension contained in $V_{\overline{P}}$ 
has the property that: among the ``free' variables, two must be equal.\/}

\medbreak
\noindent
{\bf Proof.\/} This can be seen in the structure of ${\rm CM}_{m}(\overline{d})$.
If some $d_{1i}$ are constant and the others are free but independent, the quartic terms
in the non-constant $d_{i1}$ in 
${\rm CM}_{m}(\overline{d})$ cannot vanish: their 
coefficients are non-zero by [\BERGER, Proposition 9.7.3.2] 
(the determinant, suitably normalized, is equal to the 
squared volume of the polytope formed by the points)
and the affine general position. \qed

\medbreak
\noindent

We are now ready to prove Theorem 7.1.

\medbreak
\noindent
{\bf 7.3. Proof of Theorem 7.1.} The proof is by induction on $m$. 

For $m=1$, consider $\Lambda\subset \R$ (and $\kappa\in \R_+, \kappa\ne \pm 1/2$). 
Suppose that $\Lambda$ contains two distinct points $P_1, P_2$, 
where we may as well assume $P_1=0$ and 
$P_2=R_0>0$. We consider the map sending
a third point $Q$ to the point 
$\big(d(Q, P_1), d(Q, P_2)\big)\in \R^2$. 
The image lies on the
variety $V_{\overline{P}}\subset \C^2$ 
defined by the vanishing of ${\rm CM}_1$.
The image is not weakly special 
(which is easy to see 
directly but also follows from 
Proposition 7.2), and so the number 
$\hash{\rm MSV}(W,T)$
of maximal special subvarieties 
contained in a variety $W$ with height 
up to $T$ satisfies
$$
\hash {\rm MSV}\big(V_{\overline{P}}, 
T\big) \le c(\epsilon) T^\epsilon
$$
for this family of varieties (i.e. 
uniformly in $P_0, P_1$).

Now suppose the conclusion holds for all $m'<m$ and $\Lambda\subset \R^m$. 
Suppose first that $\Lambda$ contains $m+1$ points $\overline{P}=(P_0,\ldots, P_m)$
in affine general position. Then we consider $V_{\overline{P}}$.
By Proposition 7.2, and a positive-dimensional weakly special subvariety of
$V_{\overline{P}}$ has an equality among the non-constant variables.
So any $Q$ lying on such a subvariety lies on a hyperplane determined by two
of the points $P_i \ne P_j$. There are $(m+1)m/2$ such hyperplanes, and the
induction hypothesis applies to them.
Otherwise, $Q$ is an isolated special point.

If $\Lambda$ does not contain $m+1$ points in affine general position then all the 
points in $\Lambda$ lie in some affine hyperplane, and the conclusion 
follows by induction. \ \qed

\medbreak

As noted, $V_{5, \R^3}$  does contain 
positive-dimensional non-degenerate 
weakly special subvarieties. 
It may be that all positive dimensional 
non-degenerate weakly special 
subvarieties in all dimensions 
arise in a similar way (from 
configurations that, if special, 
would be Fuglede counterexamples in 
lower dimension).

However, in $\R^2$ we can exclude 
positive-dimensional weakly special subvarieties, 
and as a consequence we show that the counterexamples to Fuglede's conjecture in $\R^2$
must be sparse. This is notionally a stronger result than the $T^\epsilon$ bound
for sets of orthogonal exponentials, though it leads only to a $T^\epsilon$ bound with
``smaller $\epsilon$''.

\medbreak

In the following, $V_{4, \C^2}$ is the variety in 
$\C^6$ defined by the vanishing
of ${\rm CM}_2$.

\medbreak
\noindent
{\bf 7.4. Proposition.\/} 
{\it The distance variety 
$V_{4, \C^2}\subset \C^6$ contains no 
positive dimensional weakly special 
subvarieties whose fixed
coordinates are all real and positive.\/}

\medbreak
\noindent
{\bf Proof.} Very elementary 
computations will be suppressed.

If $V$ contains a positive dimensional weakly special subvariety $T$,
then it contains the ``diagonal'' weakly special subvariety of $T$
in which all non-constant variables are set to be equal, and which is also
positive-dimensional. So it suffices to exclude such.

First we see that not all coordinates can be non-constant 
and equal. This
follows from the computation in \S8 below.

So some distance is fixed. Without loss of generality it is $d(P_1, P_2)$ where
$P_1=(0,0)$ and $P_2=(x_2, 0)$ is on the positive real axis.

If all other distances are non-constant and are taken to be equal, then $P_3, P_4$
lie on the bisecting (complex) line $L=\{(x_2/2, y): y\in\C\}$.
Taking a (complex) square-distance $d^2$, if $d(P_3, P_1)^2=x_3^2+y_3^2=(x_2/2)^2+y_3^2=d^2$ 
this gives two possible solutions for $P_3\in L$ and so $P_4$ must be the other one.
But $d(P_3, P_4)^2\ne d^2$ for general $d^2$.

So there must be another constant distance. Either it is a distance of $P_3$ or $P_4$
to $P_1$ or $P_2$, and without loss of generality we can assume it is $d(P_3, P_1)$,
or it is $d(P_3, P_4)$. 

In this second case, if the non-constant distances are taken equal, $P_3, P_4$
lie on $L$ (at squared distance $d^2$ to $P_1, P_2$) but are at some fixed 
distance to each other. This has only finitely many solutions for $d^2$.

We consider the first case. Now $P_3$ is at fixed (complex) distance 
$d_3$ to $P_1$, i.e.  $d_3^2=x_3^2+y_3^2$,
and we consider $P_4$ at non-constant but equal distance to $P_1, P_2, P_3$, 
so $P_4\in L$. Suppose the common ``mobile'' square-distance is $d^2$.
So $P_4=(x_2/2, y_4)$ with $y_4^2=d^2-x_2^2/4$. Then 
$$
d^2=d^2(P_4, P_3)=(x_2/2-x_3)^2+(y_4-y_3)^2
$$
and $x_3^2+y_3^2=d_3^2$ being fixed we find: the condition $d(P_3, P_2)^2=d^2$
now forces $P_3$ into a finite set.

So either $d(P_3, P_2)$ is fixed, forcing $P_4$ (equidistant to all of them) into
a finite set. Or we can assume $d(P_4, P_1)$ is also fixed and consider the three
remaining distances to be mobile but equal. 
But this is not possible.
\ \qed


\section{Transcendence conjectures and Fuglede's
conjecture}

\medskip

\noindent{\bf Schanuel's conjecture and 
non-generic Bessel zeros}

\medskip
Transcendence properties of values of the 
exponential function are
governed by Schanuel's conjecture.

\medskip
\noindent
{\bf 8.1. Schanuel's conjecture.\/} ([\LANG, p30])
Let $y_1,\ldots, y_m\in \C$. Then
$$
{\rm td}_\Q\Q(y_1,\ldots, y_m, e^{y_1},\ldots, e^{y_m}) \ge m
$$
unless $y_1,\ldots, y_m$ are linearly dependent over $\Q$.

\medbreak

We show that, for odd $m\ne \pm 1$, Schanuel's conjecture implies 
that the positive roots of $J_{m/2}$ are 
algebraically independent over $\Q$, and that this would 
confirm Fuglede's conjecture for sets of orthogonal exponentials in $\R^m$, 
for odd $m\ge 3$. (Independence fails for $J_{\pm 1/2}$.)

\medbreak
\noindent
{\bf 8.2. Proposition.\/} {\it Assume 
Schanuel's  conjecture. Let $m\ge 1$ 
be an integer. Then the positive roots 
of $J_{m+1/2}$ are algebraically
independent over $\Q$.}

\medbreak
\noindent
{\bf Proof.\/} Let $m\ge 1$. If $x$ is 
a root of $J_{m+1/2}(x)$ then we have
$$
-1=-R(x)e^{2ix}
$$
where $R$ is the rational function in 
the proof of Lemma 5.2. We know that
$R$ is non-constant and satisfies 
$R(x)R(-x)=1$. Put $y=2ix$ and define 
$S(y)=R(y/2i)^{-1}$. Then $y$
corresponding to roots satisfy $S(y)=e^y$
and  $S(y)S(-y)=1$, and we must show 
algebraic independence of roots 
$y_i$ up to $y_i=\pm y_j$.

Schanuel's conjecture immediately implies 
that any non-zero root $y_i$ is 
transcendental.

Now we suppose that any $n$ non-zero 
roots $y_i, i=1,\ldots, n$ with 
$y_i\ne\pm y_j$ for $i\ne j$ are
algebraically independent over $\Q$ 
and deduce the same for $n+1$. Let 
$y_{n+1}$ be a  non-zero root.
By Schanuel's conjecture, 
$$
{\rm td}\big(y_1,\ldots, y_{n+1}, e^{y_1},\ldots,
e^{y_{n+1}}\big)\ge n+1
$$
unless $y_1,\ldots, y_{n+1}$ are $\Q$-linearly
dependent. Each $e^{y_i}$ is algebraic over $y_i$,
and so the above reduces to
$$
{\rm td}\big(y_1,\ldots, y_{n+1}\big)\ge n+1
$$
with the same possible exception. Such a 
$\Q$ linear dependence must be of the form
$$
y_{n+1}=\sum_{j=1}^n q_jy_j
$$
for $q_j\in\Q$ not all zero, as 
$y_1,\ldots, y_n$ are algebraically 
independent by the inductive assumption.
Write $q_j=\ell_j/d$ for a positive integer 
$d$ and integers $\ell_j$.
Then
$$
S\Big(\sum_{j=1}^n q_jy_j\Big)^d=
S(y_{n+1})^d=e^{dy_{n+1}}=
\prod_{j=1}^n e^{dq_jy_j}
=\prod_{j=1}^n S(y_j)^{\ell_j}.
$$

As the $y_j$ are algebraically 
independent, considering the poles 
of the first and last displayed 
expressions, this is impossible if 
two or more $q_j$ are non-zero.
Suppose $q_j$ only is non-zero. 
Then (removing
the subscripts) we have
$$
S(qy)^d=S(y)^\ell
$$
which is impossible, considering 
degrees (as $S$ is non-constant), 
unless $\ell=\pm d$, i.e. $q=\pm 1$,
That is, unless $y_{n+1}=\pm y_j$. \ \qed

\medbreak

Concerning the algebraic independence of roots of general non-trivial 
polynomial-exponential equations and connections with Zilber's existential 
closure conjecture see [\MANTOVA, \MARKER].

\medbreak
\noindent
{\bf Algebraic independence of zeros and Fuglede's conjecture}

\medskip

Let $m\ge 1$ and $M=(m+2)(m+1)/2$. A choice of
$m+2$ points in $\R^m$ gives rise to a tuple
of distances 
$$
\overline{d}=(d_{ij})=(d(P_i, P_j), 
i\ne j)\in \R^M.
$$
Set $\ell=m+2$.
Whether a distance tuple $\overline {d}\in \R^M$
is realized by some choice of points $P_1,\ldots, P_{m+2}\in \R^m$
is governed (in part) by the vanishing of the
Cayley-Menger determinant ${\rm CM}_{m}(\overline{d})$ exhibited above.

\medbreak
\noindent
{\bf 8.3. Proposition.\/}
{\it Let $m\ge 1$ and $M=(m+2)(m+1)/2$.
Suppose that $x_i, i=1,\ldots, M$ are positive 
numbers that are algebraically independent over $\Q$.
Then there do not exist $m+2$ distinct points in $\R^m$
whose non-zero distances are among the $x_i$
(allowing repetition).}

\medbreak
\noindent{\bf Proof.\/} Without allowing repetitions
the conclusion is clear as $P={\rm CM}_{m}(\overline{d})$ is
a non-trivial polynomial relation defined over $\Q$.
So we must consider what happens if some distances are
repeated. For example if $P=X-Y$ then $P(a,b)\ne 0$
for any algebraically independent $a,b$ but $P(a,a)=0$.
So we must consider whether under some assignment
of equal variables $P$ reduces to the zero polynomial.
If so, this will be true under the assignment
of all equal variables. The conclusion follows from
some computations set out below.\ \qed

\medbreak

We need to show that, for an indeterminate $T$,
the $(\ell+1)\times (\ell+1)$ determinant
${\rm CM}_{m}((T,\ldots, T))$,
is a non-zero polynomial
in $\Z[T]$, where ${\rm CM}_{m}(\overline{T})=$
$$
\det
\begin{pmatrix}
0 & T^2 & T^2 & \ldots & T^2 & 1 \cr
T^2 & 0 & T^2 & \ldots & T^2 & 1 \cr
T^2 & T^2 & 0 & \ldots & T^2 & 1\cr
\vdots & \vdots & \vdots & \ddots & \vdots & \vdots \cr
T^2 & T^2 & T^2 & \ldots & 0 & 1\cr
1 & 1 & 1 & 1 & 1 & 0
\end{pmatrix}
=\det
\begin{pmatrix}
0 & T^2 & T^2 & \ldots & T^2 & 1 \cr
T^2 & -T^2 & 0 & \ldots & 0 & 0 \cr
T^2 & 0 & -T^2 & \ldots & 0 & 0\cr
\vdots & \vdots & \vdots & \ddots & \vdots & \vdots \cr
T^2 & 0 & 0 & \ldots & -T^2 & 0\cr
1 & 1 & 1 & 1 & 1 & 0
\end{pmatrix}
$$
where the second expression is obtained by
subtracting the first row from rows 2 to $\ell$
in the first expression.
We then expand  down the last column to get the 
$\ell \times\ell$ determinant
$$
= \epsilon
\det\begin{pmatrix}
T^2 & -T^2  & \ldots & 0  \cr
T^2 & 0  & \ldots & 0 \cr
\vdots & \vdots  & \ddots & \vdots  \cr
T^2 & 0  & \ldots & -T^2 \cr
1 & 1  & 1 & 1 
\end{pmatrix}
= \epsilon
\det
\begin{pmatrix}
T^2 & -2T^2  & \ldots & -T^2  \cr
T^2 &  -T^2 & \ldots & -T^2 \cr
\vdots & \vdots  & \ddots & \vdots  \cr
T^2 & -T^2  & \ldots & -2T^2 \cr
1 & 0  & 0 & 0
\end{pmatrix}
$$
where $\epsilon=(-1)^\ell$ and
the second expression is obtained by subtracting
the first column from all other columns.
We then expand along the bottom row to get the now 
$(\ell-1)\times(\ell-1)$ determinant
$$
=
-\det
\begin{pmatrix}
-2T^2 & -T^2 & \ldots & -T^2  \cr
  -T^2& -2T^2 & \ldots & -T^2 \cr
\vdots & \vdots & \ddots & \vdots  \cr
 -T^2 & -T^2 & \ldots & -2T^2
\end{pmatrix}
= \epsilon T^{2(\ell-1)}
\det\begin{pmatrix}
2 & 1 & \ldots & 1  \cr
  1& 2 & \ldots & 1 \cr
\vdots & \vdots & \ddots & \vdots  \cr
 1 & 1 & \ldots & 2
\end{pmatrix}
$$
$$
=(-1)^\ell\ell T^{2(\ell-1)}
$$
as the last displayed matrix has eigenvalues 
$\ell, 1,\ldots, 1$.
\qed
\medbreak

For even $m$, hence integer $m/2$, we see that
Fuglede's conjecture follows if the positive
zeros of $J_m$ are algebraically independent over
$\Q$, as he naturally expected should be the case.

\medbreak
\noindent{\bf A Bessel-Schanuel conjecture}

\medskip
There does not seem to be a Schanuel-type conjecture for the 
Bessel function stated in the literature. 
Motivated by our proof of Ax-Schanuel, we offer a formulation here.

In framing our conjecture, we have
as ``evidence''  just our Bessel-Ax-Schanuel and the 
Siegel-Shidlovskii results
[\SIEGEL, \SIEGELBOOK] and [\SHIDLOVSKII, Ch. 9, \S7].
These give the analogue of Lindemann-Weierstrass for 
suitably normalized Bessel functions of the first kind of generic rational 
non-integral order.
For integral order,  $Y_m$ is not an $E$-function, and
it seems that nothing is known about its transcendence at algebraic arguments.

In addition to the Lindemann-Weierstrass theorem, several other 
theorems and conjectures
might have served as guidance to Schanuel
in formulating his conjecture.
These include the Gelfond-Schneider theorem 
and Gelfond's theorem (independence of $\alpha^\beta$ and $\alpha^{\beta^2}$
for algebraic $\alpha\ne 0,1$ and algebraic $\beta$ of degree 3),
Schneider's conjectured generalization  [\SCHNEIDER, Problem 7]
to $\alpha^{\beta_i}$ for algebraic $\alpha\ne 0,1$ and $\Q$-linearly independent 
algebraic $\beta_i$,
Gelfond's announcement [\GELFOND] and conjecture [\SCHNEIDER, 5.12]
(and see also [\WALDSCHMIDT]),
and the widely expected algebraic independence of logarithms of multiplicatively
independent algebraic numbers (see e.g. Lang [\LANG, p31]). 

Schanuel did not have the Ax-Schanuel theorem for the 
exponential function as a guide:
indeed he conjectured it (see [\AX]), though there were 
earlier partial results of Skolem and Chabauty 
(see again [\AX]).

In addition to Bessel-Ax-Schanuel and Siegel-Shidlovskii, 
our formulation is in keeping with a conjecture of Andr\'e 
[\ANDRE, Conjecture 4.5.1] on the values of {\it holonomic 
arithmetic Nilsson Gevrey series} [\ANDRE, 3.3.1].

For $\kappa\in\Q$, but not a negative
integer, we have that
$$
K_\kappa(x)=\Gamma(\kappa+1)\big({x\over 2}\big)^{-\kappa} J_\kappa(x)
$$
is an E-function in the sense of (and as observed by) 
Siegel [\SIEGEL]. (For a negative
integer $\kappa\le -1$ the above is
undefined, but easily modified.)
We will therefore set,
for non-integer rational $\kappa$,
our normalized Bessel function of the 
first kind to be the multivalued function
$$
J_\kappa^*(x)=\Gamma(\kappa+1)J_\kappa(x).
$$
The rational power of $x$ will not 
affect the transcendence statement,
while $J^*_\kappa(x)$ remains a solution 
of the Bessel d.e. of order $\kappa$.

For non-negative integer $\kappa$ we set $J_\kappa^*=J_\kappa$.

Then for non-integral rational $\kappa$ 
we have a normalized fundamental 
matrix $W_\kappa^*$ formed using
$J_\kappa^*$ and $Y^*_\kappa=J_{-\kappa}^*$
(or any linear combination of 
$J^*_\kappa$ and $J^*_{-\kappa}$ over 
$\overline{\Q}$
that is independent of $J^*_\kappa$); 
with the given choice we have 
$$
{\rm Wr}(J^*_\kappa, 
Y^*_\kappa)=-2\kappa/x.
$$

For integer order $k$, restricting to 
non-negative $k$ for simplicity, we 
modify $Y_k$ as follows. For $k=0$ we 
have ([\BATEMAN, 7.24, eqn. (33)])
$$
\pi Y_0(x)= 2J_0(x)\big(\log (x/2) 
+\gamma\big) - 
2 \sum _{j=0}^\infty {(-1)^j\over j! j!} 
\big({x\over 2}\big)^{2j} H_j
$$
where $H_j=1+1/2+\ldots 1/j$ is the 
harmonic number, $H_0=0$, and $\gamma$ is 
Euler's constant.
Accordingly define
$$
Y_0^*(x)= J_0(x) \log x +A_0(x), \quad 
A_0(x)= - \sum _{j=0}^\infty {(-1)^j\over 
j! j!} \big({x\over 2}\big)^{2j} H_j.
$$
We note that $A_0$ is an E-function, 
and so $Y_0$ is a holonomic arithmetic 
Nilsson Gevrey series as defined in 
[\ANDRE]. 

For other (positive) integer orders 
we proceed similarly. We have
$$
\pi Y_n(x)=
2\Big(\gamma+\log(x/2)\Big)J_n(x)
$$
$$
-\sum_{m=0}^{n-1} {(n-m-1)!\over m!} \Big({x\over 2}\Big)^{2m-n}
-\sum_{m=0}^\infty
{(-1)^m(H_{m+n}+H_m)\over m!(n+m)!}
\Big({x\over 2}\Big)^{n+2m}.
$$
(see e.g. [\BATEMAN, 7.24, eqn. (32)]).
So we set
$$
Y_n^*(x)=(\pi/2)Y_n(x)-
(\gamma-\log 2)J_n(x).
$$
which is (but for the Laurent 
polynomial) a holonomic arithmetic 
Nilsson Gevrey series. 
Then, for non-negative integer $\kappa$,
we have
$$
{\rm Wr}(J^*_k, Y^*_k)=1/x.
$$

For a negative integer $\kappa$,
the d.e. depending only on $\kappa^2$,
we will take
$$
J_\kappa^*=J_\kappa = 
(-1)^\kappa J_{-\kappa}, 
\quad Y_\kappa^*=(-1)^\kappa Y_{-\kappa}^*.
$$

For $\kappa\in\Q$,
let $W_\kappa^*$ be the fundamental matrix
formed using $J^*_\kappa, Y^*_\kappa$.
Since the Wronskian is an algebraic 
function of x over $\overline{\Q}$, 
the maximal possible transcendence 
degree of $\Q\big(x, W^*_\kappa(x)\big)$ 
over $\Q(x)$ is 3. Here and in the 
following, any determination of the value of
a multivalued function is permitted.

\medbreak
\noindent
{\bf 8.4. Conjecture.\/}
Let $\kappa\in \Q\backslash (\Z+1/2)$.
Let $y_1,\ldots, y_m \in \C$ and 
$x_1,\ldots, x_n\in\C^\times $. Then
$$
{\rm td}_{\Q}\Q\Big(y_1,\ldots, y_m, 
e^{y_1},\ldots, e^{y_m},
x_1,\ldots, x_n,
W^*_\kappa(x_1),\ldots, W^*_\kappa(x_n)\Big)\ge m+3n,
$$
unless either

(A) $y_1, \ldots, y_m$ are $\Q$-linearly dependent, or

(B) $x_i=\pm x_j$ for some $i\ne j$.

\medbreak

Note that, while we normalized to remove $\pi$ from the expansion of the $Y_n$, as 
Conjecture 8.4 includes Schanuel's Conjecture (the case $n=0$), 
algebraic independence of Bessel arguments and values over $\Q(\pi)$ is implied. 
More generally, a transcendence degree of $3n$ for 
the Bessel arguments and values is implied (assuming 
$x_1,\ldots, x_n$ are $\pm$ distinct)
over the (exponential) field $E$ of ``exponential algebraic numbers''; see [\KIRBY]. For a related conjecture (with a related provenance)
see 
[\MANGOUBIROSENBLATT].
One can make a corresponding formulation involving normalized
Bessel functions of several different generic rational orders $\kappa_1,\ldots, \kappa_\ell$, 
provided $\kappa_i \pm \kappa_j \not\in \Z$ for $i\ne j$. Evaluating at 
$y_j\in \C, x_i\in\C^\times$ as above, and with the same caveats, 
one expects transcendence degree at least $m+3n\ell$.

\bigbreak

\centerline{\bf References}

\begin{enumerate}

\item[\ANDRE.] Y. Andr\'e, Arithmetic Gevrey series and transcendence. A survey,
{\it J. Th. Nombres Bordeaux\/} {\bf 15} (2003), 1--10.


\item[\AWH.] G. B. Arfken, H. J. Weber, and F. E. Harris, 
{\it Mathematical methods for physicists: a comprehensive guide,\/} 
Academic Press, 2011.


\item[\AX.] J. Ax, On Schanuel’s conjectures, 
{\it Annals Math. (2)\/} {\bf 93} (1971), 252–268.


\item[\BALSER.] W. Balser,  Formal power series and linear systems of meromorphic 
ordinary differential equations, {\it Springer Science \& Business Media\/}, 2008


\item[\BATEMAN.] Staff of the Bateman Manuscript Project,
{\it Higher transcendental functions,}
McGraw-Hill, New York -- Toronto -- London, 1953.


\item[\BEURLING.] A. Beurling, 
{\it Quasianalyticity. Mittag-Leffler Lectures on Complex Analysis (1977-1978), 
Collected Works, Vol I,\/} Birkh\"auser, Boston, 1989.


\item[\BERGER.] M. Berger, {\it Geometry I\/}, 
Corrected fourth printing, Springer-Verlag, 
Berlin--Heidelberg, 2009.






\item[\BKPLMM.] G. Binyamini, A. Kiro, and J. Pila,
Laguerre-Manin-Mumford, in preparation.


\item[\BNZ.]
G. Binyamini, D. Novikov, and B. Zak, 
Wilkie's conjecture for Pfaffian structures
{\it Ann. of Math. (2)} 
{\bf 199} (2024), 795–821.


\item[\BCFNDAS.]
D. Bl\'azquez-Sanz, G. Casale, J. Freitag,  and J. Nagloo,
A  differential approach to the 
Ax-Schanuel, I, 
arXiv:2102.03384.




\item[\DLMF.] {\it NIST Digital Library of 
Mathematical Functions,\/}  Release 1.2.4 
of 2025-03-15. F. W. J. Olver, 
A. B. Olde Daalhuis, D. W. Lozier, 
B. I. Schneider, R. F. Boisvert, C. W. Clark, 
B. R. Miller, B. V. Saunders, H. S. Cohl,
and M. A. McClain, editors. Available at
https://dlmf.nist.gov


\item[\DRIESSPEISSEGGER.]
L. van den Dries and P. Speissegger, 
The field of reals with multisummable series 
and the exponential function,
{\it Proc. London Math. Soc. (3)\/}  {\bf 81} 
(2000), 513–-565.


%


\item[\DYNKIN.] E. M. Dyn'kin, The pseudoanalytic extension, 
{\it Journal d'Analyse Math\'ematique\/} {\bf 60} (1993), 45--70.


\item[\FUGLEDE.] B. Fuglede,
Commuting self-adjoint partial differential 
operators and a group theoretic problem, 
{\it J. Functional Anal.\/} {\bf 16} (1974), 
101–121.


\item[\FUGLEDEBALL.] B. Fuglede,
Orthogonal exponentials on the ball, 
{\it Expo. Math.} {\bf 19} (2001), 267–272.


\item[\GELFOND.] A. O. Gelfond,
Sur quelques r\'esultats nouveaux dans la th\'eorie des nombres transcendants,
{\it C. R. Math. Acad. Sci. Paris\/} {\bf 199} (1934), 
259, available at: http://gallica.bnf.fr/ark:/12148/  
bpt6k3151h.image.r=gelfond.langEN.f259.pagination




\item[\IK.]
A. Iosevich and M. N. Kolountzakis,
Size of orthogonal sets of exponentials 
for the disk,
{\it Rev. Mat. Iberoam.\/} {\bf 29} (2013), 
739–747.


\item[\JK.]
R. Jaoui and J. Kirby,
Integration in finite terms and  
exponentially algebraic functions,
arXiv:2510.26248.

\item[\KIRBY.] J. Kirby, 
Exponential algebraicity in 
exponential fields,
{\it Bull. Lond. Math. Soc.\/} 
{\bf 42} (2010), 879–890.



\item[\KOLCHIN.] E. Kolchin,
Algebraic groups and algebraic dependence, 
{\it Amer. J. Math.\/} {\bf 90}
(1968), 1151--1164.


\item[\LANG.] S. Lang, 
{\it Introduction to transcendental numbers,\/}
Addison-Wesley Publishing Company, 
Reading, Massachusetts -- Palo Alto -- London -- Don Mills, Ontario, 1966.


\item[\MANGOUBIROSENBLATT.] D. Mangoubi and D. Rosenblatt,
On multiplicity bounds for eigenvalues of the clamped round plate,
arXiv:2411.07564.


\item[\MANTOVA.]
V. Mantova,
Polynomial-exponential equations and Zilber's conjecture,
With an appendix by Mantova and U. Zannier,
{\it Bull. Lond. Math. Soc.\/} {\bf 48} (2016), 309–320.


\item[\MARKER.]
D. Marker, 
A remark on Zilber's pseudoexponentiation.
{\it J. Symbolic Logic\/} {\bf 71} (2006), 791–798. 




\item[\OLVER.] F. Olver, {\it Asymptotics and 
special 
functions,\/} AK Peters/CRC Press, 1997.


\item[\PW.] J. Pila and A. J. Wilkie,
The rational points of a definable set,
{\it Duke Math. J.} {\bf  133} (2006), 591–616.


\item[\PZ.] J. Pila and U. Zannier,
Rational points in periodic analytic sets and the 
Manin-Mumford conjecture,
{\it Atti Accad. Naz. Lincei Rend. 
Lincei Mat. Appl.\/}
{\bf 19} (2008), 149–162.

\item[\PUTSINGER.]
M. van der Put and M. Singer,
{\it Galois theory of linear differential 
equations,\/}
Grundlehren Math. Wiss., {\bf 328}, 
Springer-Verlag, 
Berlin, 2003.


\item[\RAYNAUD.] M. Raynaud, 
Sous-vari\'et\'es d'une vari\'et\'e ab\'elienne et
points de torsion, in {\it Arithmetic and Geometry, 
Vol. I}, 
Prog. Math {\bf 35}, M. Artin and J. Tate, eds,
Birkh\"auser, Boston, 1983, 327--352.


\item[\SCHNEIDER.] T. Schneider, 
{\it  Einf\"uhrung in die transzendenten Zahlen},
Springer-Verlag, Berlin-G\"ottingen-Heidelberg, 1957.


\item[\SHIDLOVSKII.] A. B. Shidlovskii, {\it Transcendental numbers,\/} 
Translated from the Russian by N. Koblitz, Walter de Gruyter, Berlin-New York, 1989.


\item[\SIEGEL.] C. L. Siegel, 
\"Uber einige Anwendungen diophantischer
Approximationen, {\it Abh. der preussischen Akad. 
der Wissenschaften, phys-math. Klasse,} 
Jahrgang 1929, Nr. 1 (1930), 1--70. 
English translation by C. Fuchs in:
{\it On some applications of 
Diophantine approximations\/},
Quad./Monogr. {\bf 2}, U. Zannier, editor,
Edizioni della Normale,  Pisa, 2014, 1–80.


\item[\SIEGELBOOK.] C. L. Siegel,
{\it Transcendental numbers,\/}
Annals of Math. Studies {\bf 16}, 
Princeton University Press, 1949.




\item[\WALDSCHMIDT.]
M. Waldschmidt,
Conjectures and open problems, 
Lecture 5, Arizona Winter School, Tucson, Arizona, \ 2008, 
slides, at: https: //swc-math.github.io/aws/2008/08WaldschmidtSlides5.pdf


\item[\WASOW.] W. Wasow,
{\it Asymptotic expansions for ordinary 
differential equations,} 
Courier Dover Publications, 2018.


\item[\ZAKHAROV.] D. Zakharov,
On sets of orthogonal exponentials on the disk,
arXiv: 2405.14063.

\end{enumerate}

\bigskip
\noindent
Dept. of Mathematics,
Weizmann Institute of Science, Rehovot, Israel

\medskip
\noindent
and (GB and JP)

\medskip
\noindent
School of Mathematics, Institute for Advanced Study, 
Princeton, USA

\bigskip

\noindent
gal.binyamini@weizmann.ac.il

\smallskip

\noindent
avner-ephraiem.kiro@weizmann.ac.il 

\smallskip

\noindent
gady.kozma@weizmann.ac.il

\smallskip

\noindent
jonathan.pila@maths.ox.ac.uk


\end{document}